\documentclass[10pt]{article}
\usepackage{lmodern}
\usepackage{amsmath}
\usepackage[T1]{fontenc}
\usepackage[utf8]{inputenc}
\usepackage{authblk}
\usepackage{amsfonts}
\usepackage{graphicx}
\usepackage{tkz-euclide}
\tikzset{elegant/.style={smooth,thick,samples=50,cyan}}
\tikzset{eaxis/.style={->,>=stealth}}
\usepackage{color,tikz}
\usetikzlibrary{decorations.pathreplacing,calc}
\usepackage{rotating}
\usepackage{amssymb}
\usepackage[english]{babel}
\usepackage{color}
\usepackage{amsthm}
\usepackage{graphicx}
\usepackage{mathrsfs}
\usepackage{makecell}
\usepackage{microtype}
\usepackage{mathscinet}
\usepackage{array}
\usepackage{multirow}
\usepackage{enumerate}
\usepackage[cal=boondoxo,bb=ams]{mathalfa}
\usepackage{hyperref}
\hypersetup{hidelinks}
\usepackage{booktabs}
\usepackage{arydshln}

\newtheorem{theorem}{Theorem}[section]
\newtheorem{prop}{Proposition}[section]

\newtheorem{remark}{Remark}[section]

\newcommand{\ml}{\mathcal}
\newcommand{\mb}{\mathbb}

\DeclareMathOperator{\intt}{int}
\DeclareMathOperator{\extt}{ext}
\DeclareMathOperator{\bdd}{bdd}
\DeclareMathOperator{\diag}{diag}

\title{Some asymptotic profiles for the viscous Moore-Gibson-Thompson equation in the $L^q$ framework}
\author[1]{Wenhui Chen\thanks{Wenhui Chen (wenhui.chen.math@gmail.com)}}
\affil[1]{School of Mathematics and Information Science, Guangzhou University, 510006 Guangzhou, China}
\author[2]{Junying Gong \thanks{Junying Gong (2657683160@qq.com)}}
\affil[2]{School of Medical Information Engineering, Guangdong Pharmaceutical University, 510006 Guangzhou, China}
\date{}

\begin{document}

		\maketitle
		\begin{abstract}
			\medskip
	This manuscript studies some qualitative properties of solutions to the Cauchy problem for the viscous Moore-Gibson-Thompson (MGT) equation. For one thing, by applying the WKB analysis and diagonalization procedure, we derive some $L^p-L^q$ decay estimates and the large time asymptotic profile for a suitable energy term, which provides a new way to treat higher order MGT-type coupled systems. For another, we obtain the global (in time) singular limits in the $L^q$ framework and the higher order asymptotic profile with respect to small thermal relaxation via the multi-scale analysis and the Fourier analysis. Especially, provided the incompatible initial condition between the viscous MGT equation and the strongly damped wave equation, the formation of initial layer is rigorously justified.
			\\
			
			
			\noindent\textbf{AMS Classification (2020)}  35L30, 35B40, 35B25
		\end{abstract}
\fontsize{12}{15}
\selectfont

\section{Introduction}\label{Section-Introduction}\setcounter{equation}{0}
\hspace{5mm}In the present paper, we consider the Cauchy problem for the viscous Moore-Gibson-Thompson (MGT) equation, arising from a linearized model for acoustic waves propagation in viscous thermally relaxing fluids, namely,
\begin{align}\label{Eq_MGT}
	\begin{cases}
		\tau\psi_{ttt}+\psi_{tt}-\Delta \psi-(\delta+\tau)\Delta\psi_t=0,&x\in\mb{R}^n,\ t>0,\\
		\psi(0,x)=\psi_0(x),\ \psi_t(0,x)=\psi_1(x),\ \psi_{tt}(0,x)=\psi_2(x),&x\in\mb{R}^n,
	\end{cases}
\end{align}
with the thermal relaxation $\tau>0$ and the diffusivity of sound $\delta>0$ containing viscous coefficients of a given thermoviscous flow (cf. Remark \ref{Rem-viscosity}), where the unknown function $\psi=\psi(t,x)\in\mb{R}$ is referred to the acoustic velocity potential in the classical theory of acoustic waves. Our main purpose is to study asymptotic profiles for the viscous MGT equation \eqref{Eq_MGT} in depth with large time $t\gg1$ or small thermal relaxation $0<\tau\ll 1$. In particular, we formulate a new singular (initial) layer with respect to small thermal relaxation under the incompatibility of initial conditions between the cases $\tau>0$ and $\tau=0$.

It is widely known that to characterize the propagation of sound in viscous thermally relaxing fluids, the classical study of nonlinear acoustics with second sound phenomenon (cf. the pioneering effort \cite{Blackstock-1963,Hamilton-Blackstock-1998,Jordan-2014}) considers an approximated model of the fully compressible Navier-Stokes-Cattaneo system in irrotational flows. To be specific, by using the Lighthill scheme of approximations to retain the terms of first and second orders with small perturbations around the equilibrium state, the following well-studied Jordan-MGT equations (see \cite{Kaltenbacher-Lasiecka-Pos-2012,Kaltenbacher-Nikolic-2019,Racke-Said-2020,B-L-2020,Said-2021,Said-2021-Bes,Kaltenbacher-Nikolic-2021,Kaltenbacher-Nikolic-2022,Chen-Takeda=2023} and references therein) occurs:
\begin{align}\label{JMGT}
	\tau\psi_{ttt}+\psi_{tt}-\Delta\psi-(\delta+\tau )\Delta\psi_t=\partial_t\left(\frac{B}{2A}|\psi_t|^2+|\nabla\psi|^2\right),
\end{align}
with the coefficients of nonlinearity fulfilling $B/A>0$. It's worth noting that $\delta>0$ is always referred to the viscous case due to its viscous dissipation partly deriving from the Navier-Stokes equations. The corresponding linearization of the viscous Jordan-MGT equation \eqref{JMGT} with the vanishing right-hand side, i.e. the viscous MGT equation (see \cite{Kaltenbacher-Lasiecka-Marchand-2011,Marchand-McDevitt-Triggiani-2012,Conejero-Lizama-Rodenas-2015,Dell-Pata-2017,Pellicer-Said-Houari=2019,B-C-Lasiecka=2020,Chen-Ikehta=2021} and references therein), is given by 
\begin{align}\label{MGT} 
	\tau\psi_{ttt}+\psi_{tt}-\Delta\psi-(\delta+\tau )\Delta\psi_t=0,
\end{align}
which has been proposed in the early literature of F.K. Moore, W.E. Gibson \cite{MooreGibson1960} in 1960 and P.A. Thompson \cite{Thompson1972} in 1972. Consequently,  the third order (in time) strictly hyperbolic equation \eqref{MGT} is called the viscous MGT equation. These scalar equations and their related mathematical models  have been extensively considered in medical and industrial applications of high-intensity
ultra sound, for instance, medical imaging and therapy, ultrasound cleaning and welding (see \cite{Abramov-1999,Dreyer-Krauss-Bauer-Ried-2000,Kaltenbacher-Landes-Hoffelner-Simkovics-2002} and references given therein).

To deeply understand the underlying physical phenomena in a certain condition, it is important to study qualitative properties of solutions to the linearized model. Let us focus on the Cauchy problem for the viscous MGT equation \eqref{Eq_MGT}, which has been firstly studied by the paper \cite{Pellicer-Said-Houari=2019}. The authors of \cite{Pellicer-Said-Houari=2019} employed energy methods in the Fourier space combined with suitable Lyapunov functionals to derive energy estimates, and eigenvalues expansions to investigate some $L^2$ estimates for the solution itself. Later, by applying the explicit representation of solutions and the Fourier analysis, \cite{Chen-Ikehta=2021} obtained optimal growth ($n=1,2$) and decay ($n\geqslant3$) estimates of solutions for large time in the $L^2$ framework. Moreover, they also got local (in time) singular limits as thermal relaxation tending to zero when $n\geqslant 3$ with the aid of classical energy methods and Hardy's inequality. In the recent paper \cite{Chen-Takeda=2023}, the authors made use of the asymptotic analysis to capture the optimal leading terms and second order large time profiles. Summarizing the used tools on large time behavior for the Cauchy problem \eqref{Eq_MGT}, they are mainly based on energy methods or explicit representations of solutions. Nevertheless, several kinds of higher order MGT-type coupled systems, arising in acoustics, mechanics, heat transfer and heat conduction equations (cf. \cite{Baz-Fer-Qui=2021,Baz-Fer-Mag-Qui=2022,Baz-Fer-Qui=2023,Baz-Fer-Liv-Qui=2024}), are challenging to be studied in the $L^q$ framework by the method mentioned above (see Remark \ref{Rem-Diag-Large-Class} for the reasons in detail). Hence, it is necessary to develop a unified treatment of MGT-type models. The viscous MGT equation \eqref{Eq_MGT}, as a fundamental MGT-type model, will be treated by this manuscript in the $L^q$ framework, which provides a new way to study a large class of MGT-type coupled systems.

Importantly, one may notice that the highest time-derivative of solution and the third initial condition in the Cauchy problem \eqref{Eq_MGT} will lose as $\tau=0$. If the incompatible initial condition does not hold additionally, it may generate a singular layer, i.e. an initial layer as $\tau\downarrow0$. As is known to all, rigorous mathematical studies for singular limits and singular layers are still challenging in fluid PDEs, for example, the Sobolev stability and instability for the Prandtl boundary layer equation \cite{Gerard-Dormy=2010,Alexandre-Wang-Xu-Yang=2015}. Therefore, it seems quite interesting to study singular limits and asymptotic profiles of solutions in the initial layer to the viscous MGT equation \eqref{Eq_MGT} with small thermal relaxation $0<\tau\ll 1$. Note that the singular limits for the Cauchy problem \eqref{Eq_MGT} as $\tau\downarrow 0$ has been completed (cf. \cite{B-C-Lasiecka=2020,Chen-Ikehta=2021,Chen-Takeda=2023} with the linear rate $\tau$ of convergence). In other words, the solution $\psi=\psi(t,x)$ for the viscous MGT equation \eqref{Eq_MGT} converges to the solution $\psi^0=\psi^0(t,x)$ for the strongly damped wave equation (or the so-called viscoelastic damped wave equation) as follows:
\begin{align}\label{Eq_VEDW}
	\begin{cases}
		\psi_{tt}^0-\Delta \psi^0-\delta\Delta\psi_t^0=0,&x\in\mb{R}^n,\ t>0,\\
		\psi^0(0,x)=\psi_0(x),\ \psi_t^0(0,x)=\psi_1(x),&x\in\mb{R}^n.
	\end{cases}
\end{align}
 It seems that higher order global (in time) singular limits and the formulation of initial layer for the viscous MGT equation \eqref{Eq_MGT} are still unknown. We will partly answer these questions later.

Our contribution of this manuscript is twofold. For one thing, by applying the refined WKB analysis and multi-step diagonalization procedure (developed in \cite{Yagdjian=1997,Reissig-Wang=2005,Yang-Wang=2006,Jachmann=2008,Reissig=2016}) in Section \ref{Section-Lp-Lq}, we obtain the sharp $L^p-L^q$ decay estimates for the new energy term $\Psi=\Psi(t,x)\in\mb{R}^3$ which is defined via
\begin{align}\label{Energy-Term-Psi}
	\Psi:=\big(\psi_t+ |D|\psi,\psi_t-|D|\psi,\tau\psi_{tt}-(\delta+\tau)\Delta\psi\big)^{\mathrm{T}}.
\end{align}
The pseudo-differential operator $|D|$ carries its symbol $|\xi|$. Furthermore, the large time profile of $\Psi$ is related to $\Phi=\Phi(t,x)$ for the first order (in time) system
\begin{align}\label{Ref-System}
\Phi_t+\frac{1}{\tau}\diag(0,0,1)\Phi+\diag(i,-i,0)|D|\Phi-\frac{\delta}{2}\diag(1,1,-2)\Delta\Phi=0
\end{align}
with suitable initial data $\Phi_0=\Phi_0(x)$ that will be chosen later. This approach, particularly the new energy term $\Psi$, provides a new way to study a large class of MGT-type coupled systems. For another, by using the multi-scale analysis and the Fourier analysis in Section \ref{Section-Singular-Limit}, we rigorously justify the following global (in time) convergence result:
\begin{align*}
\psi\to\psi^0+\tau\psi^{I,1}+\tau^2\left(\psi^{I,2}+\psi^{\mathrm{Lay}}\right)\ \ \mbox{in}\ \ L^{\infty}([0,+\infty),L^q)
\end{align*}
when $\tau\downarrow0$ with the rate $\tau^3$ of convergence, where the higher order profiles $\psi^{I,1}$ and $\psi^{I,2}$ are the solutions to the inhomogeneous strongly damped wave equations \eqref{Kuznetsov_I_1} and \eqref{Kuznetsov_I_2}, respectively. Moreover, the crucial layer function $\psi^{\mathrm{Lay}}$ with a fast change factor as $\tau\downarrow0$ is defined in \eqref{fun-lay}, which is non-trivial under the incompatible initial condition $\psi_2\neq\Delta\psi_0+\delta\Delta\psi_1$. This result not only shows the small thermal relaxation profile of $\psi$, but also formulates the new initial layer rigorously in the $L^{\infty}([0,+\infty),L^q)$ norm. Note that some comments for the inviscid MGT equation with $\delta=0$ are addressed in Section \ref{Sec-Final-Rem}.

\section{Main results}\label{Section-Main-Results}\setcounter{equation}{0}
\subsection{Large time profile in the $L^q$ framework}
\hspace{5mm}Our first result shows $L^q$ the well-posedness, $L^p-L^q$ decay estimates and large time asymptotic profile of the energy term $\Psi$ defined in \eqref{Energy-Term-Psi} to the viscous MGT equation \eqref{Eq_MGT}. Note that this energy term consists of some derivatives of the acoustic velocity potential, that is $|D|\psi$, $\Delta\psi$, $\psi_t$ and $\psi_{tt}$, whose motivation of its introduction will be addressed in \eqref{Motivation}. Before stating the main theorem, let us take some auxiliary matrices as follows:
\begin{align}
	T_{1,s}&:=\left(\begin{array}{ccc}
		1 & 0 & -1\\
		0 & 1 & -1\\  
		0 & 0 & 1
	\end{array}\right),\ \ N_{2,s}:=\frac{1}{2}\left(\begin{array}{ccc}
	0 & 0 & 2\tau\\
	0 & 0 & -2\tau\\  
	\delta & -\delta & 0
\end{array}\right),\ \ T_{3,s}:=\left(\begin{array}{ccc}
i & -i & 0\\
1 & 1 & 0\\  
0 & 0 & 1
\end{array}\right),\label{T1,N2,T3}\\
N_{3\frac{1}{2},s}&:=\frac{1}{2}\left(\begin{array}{ccc}
	0 & 0 & (1-i)\tau^2\\
	0 & 0 & (1+i)\tau^2\\  
	-(1+i)\tau\delta & -(1-i)\tau\delta & 0
\end{array}\right), \ \ N_{4,s}=\frac{1}{4}\left(\begin{array}{ccc}
0 & -\delta  & 0\\
-\delta  & 0 & 0\\
0 & 0 & 0
\end{array}\right),\label{N3-1/2,N4}
\end{align}
whose constructions will be shown in Subsection \ref{Sub-Section-Diag-Small}. Moreover, we denote the identity matrix in three dimensions by $I_{3}:=\diag(1,1,1)$, and the diagonal matrix with three exponential elements by
\begin{align*}
\diag\left(\mathrm{e}^{\lambda_j(|\xi|)t}\right)_{j=1}^3:=\diag\left(\mathrm{e}^{\lambda_1(|\xi|)t},\mathrm{e}^{\lambda_2(|\xi|)t},\mathrm{e}^{\lambda_3(|\xi|)t}\right).
\end{align*}
Let us define the pseudo-differential operator in a matrix sense as follows:
\begin{align*}
\ml{T}(|D|):=T_{1,s}(I_{3}+N_{2,s}|D|)T_{3,s}(I_{3}+N_{3\frac{1}{2}}|D|^2).
\end{align*}
Later,  $f\lesssim g$ means $f\leqslant Cg$ with a positive constant $C$ changing from line to line.
\begin{theorem}\label{Thm-Lp-Lq}
Suppose that the initial data $\Psi_0\in(H^{s+M_{p,q,n}}_p)^3$ with $s\geqslant 0$, $1\leqslant p\leqslant 2\leqslant q\leqslant+\infty$, and
\begin{align}\label{Mpq}
M_{p,q,n}\begin{cases}
	>n(\frac{1}{p}-\frac{1}{q})&\mbox{when}\ \ p\neq q,\\
	=0&\mbox{when}\ \ p=q.
\end{cases}
\end{align} Then, there is a unique determined Sobolev solution in the sense of the suitable energy term to the viscous MGT equation \eqref{Eq_MGT} such that
\begin{align*}
\Psi\in\big(\ml{C}([0,+\infty),H^s_q)\big)^3.
\end{align*}
It satisfies the $L^p-L^q$ decay estimates
\begin{align*}
\|\Psi(t,\cdot)\|_{(\dot{H}^s_q)^3}\lesssim (1+t)^{-\frac{s}{2}-\frac{n}{2}(\frac{1}{p}-\frac{1}{q})}\|\Psi_0\|_{(H^{s+M_{p,q,n}}_p)^3},
\end{align*}
 as well as the refined estimates
\begin{align*}
\left\|\Psi(t,\cdot)-\ml{T}(|D|)\Phi(t,\cdot)\right\|_{(\dot{H}^s_q)^3}\lesssim(1+t)^{-\frac{1}{2}-\frac{s}{2}-\frac{n}{2}(\frac{1}{p}-\frac{1}{q})}\|\Psi_0\|_{(H^{s+M_{p,q,n}}_p)^3},
\end{align*}
in which $\Phi$ is the solution to the reference system \eqref{Ref-System} with its initial data $\Phi_0:=\ml{T}^{-1}(|D|)\Psi_0$.
\end{theorem}
\begin{remark}
By subtracting the function $\ml{T}(|D|)\Phi(t,\cdot)$ in the $(\dot{H}^s_q)^3$ norm, the decay rate of $\Psi(t,\cdot)$ can be improved $(1+t)^{-\frac{1}{2}}$ so that
\begin{align*}
\lim\limits_{t\to+\infty}t^{\frac{s}{2}+\frac{n}{2}(\frac{1}{p}-\frac{1}{q})}\left\|\Psi(t,\cdot)-\ml{T}(|D|)\Phi(t,\cdot)\right\|_{(\dot{H}^s_q)^3}=0.
\end{align*}
 Thus, we may explain the  large time profile of the energy term $\Psi$ by the diffusion waves function $\Phi$. Motivated by its large time profile $\Phi$ satisfying \eqref{Ref-System}, particularly the viscous term $-\delta\Delta\Phi$, the decay property of $\Psi$ is generated by the viscous effect with $\delta>0$. This phenomenon coincides with the nature of viscosity in the viscous MGT equation and the previous researches \cite{Pellicer-Said-Houari=2019,Chen-Ikehta=2021,Chen-Takeda=2023}.
\end{remark}

\begin{remark}\label{Rem-Prop-sing}
The property of propagation of $H^s$-singularities (i.e. mild singularities) for solutions to the viscous MGT equation \eqref{Eq_MGT} could be of interest. Note that a function $v=v(x)$ belongs to $H^s_{\mathrm{loc}}(\{x_0\})$ if there exists a neighborhood $\ml{U}_{\epsilon_0}(x_0)=\{x\in\mb{R}^n:|x-x_0|<\epsilon_0\}$ such that $\langle\xi\rangle^s\ml{F}(\chi v)\in L^2$ for all $\chi\in\ml{C}_0^{\infty}$ with $\mathrm{supp}\,\chi\subset\ml{U}_{\epsilon_0}(x_0)$. Moreover, a function $v\in H^s_{\mathrm{loc}}(\Omega)$ with $\Omega\subset\mb{R}^n$ provided that $\chi v\in H^s$ for all $\chi\in\ml{C}_0^{\infty}(\Omega)$.

Let us turn back to the propagation of $H^s$-singularities for our problem \eqref{Eq_MGT}. Suppose that $\psi_1\pm|D|\psi_0$, $\tau\psi_2-(\delta+\tau)\Delta\psi_0$ belongs to $H^s$ but not $H^{s+1}_{\mathrm{loc}}(x_0)$ for a given point $x_0\in\mb{R}^n$. Then,
\begin{align*}
\psi_t(t,\cdot),|D|\psi(t,\cdot)\not\in H^{s+1}_{\mathrm{loc}}(\{x_0\pm c_0te_0\}) \ \ \mbox{for all}\ \ t>0
\end{align*}
with the propagation speed $c_0:=\sqrt{\frac{\delta+\tau}{\tau}}>0$, where $e_0$ is an arbitrary unite vector in $\mb{R}^n$. This propagation speed coincides with the one in \cite[Remark 3.4]{Chen-Ikehta=2021} by classical energy methods. Its proof is based on Proposition \ref{Prop-large-frenquencies} and \cite[Theorem 2.26]{Jachmann=2008}.
\end{remark}


\begin{remark}\label{Rem-Diag-Large-Class}
The mathematical models based on the viscous MGT equation \eqref{Eq_MGT} have deserved much interests in the last years, including the context of acoustics, mechanics, heat transfer and heat conduction equations in the thermoelastic theory (see, for example, \cite{Baz-Fer-Qui=2021,Baz-Fer-Mag-Qui=2022,Baz-Fer-Qui=2023,Baz-Fer-Liv-Qui=2024} and references given therein). To study large time asymptotic behavior of solutions for these models in depth, nevertheless, the approaches of \cite{Chen-Ikehta=2021,Chen-Takeda=2023} are quite difficult to apply because tedious and complicated computations of characteristic roots as well as explicit solution's formulas for higher order PDEs. Thus, one needs another approach to treat higher order MGT-type  systems based on the viscous MGT equation.

By introducing a new energy term $\Psi\in\mb{R}^3$ for the viscous MGT equation \eqref{Eq_MGT}, we derive some large time qualitative properties of solutions in the $L^q$ framework. Since the viscous MGT equation \eqref{Eq_MGT} is the fundamental model of a large class of MGT-type models, we believe that our approach can be widely applied in these MGT-type models rather than the explicit analysis in \cite{Chen-Ikehta=2021,Chen-Takeda=2023}. For example, one may define a suitable energy term related to \eqref{Energy-Term-Psi} to study the Cauchy problem for the MGT-Fourier model, which describes the vibrations of a viscoelastic heat conductor obeying the Fourier thermal law and governed by the standard linear solid model (cf. \cite{Alves-Buriol-Ferr-Mun=2013,Conti-Liv-Pata=2021,Del-Pata=2023} and references therein), as follows:
\begin{align*}
\begin{cases}
\psi_{ttt}+\psi_{tt}-\Delta\psi-\beta\Delta\psi_t=-\eta\Delta\theta,\\
\theta_t-\Delta\theta=\eta\Delta\psi_{tt}+\eta\Delta\psi_t,
\end{cases}
\end{align*}
with the viscous parameter $\beta>0$ and the coupling constant $\eta\neq0$.
\end{remark}

\subsection{Small thermal relaxation profile in the $L^q$ framework}
\hspace{5mm}We now turn to singular limits for the Cauchy problem \eqref{Eq_MGT} with respect to small thermal relaxation $0<\tau\ll 1$. Let us recall the formal limit model as $\tau=0$, that is the strongly damped wave equation \eqref{Eq_VEDW} with its solution $\psi^0$. Moreover, we introduce two functions $\psi^{I,1}=\psi^{I,1}(t,x)$ and $\psi^{I,2}=\psi^{I,2}(t,x)$ for the auxiliary Cauchy problems, respectively, as follows:
\begin{align}\label{Kuznetsov_I_1}
	\begin{cases}
		\psi_{tt}^{I,1}-\Delta\psi^{I,1}-\delta\Delta\psi_t^{I,1}=\Delta\psi_t^0-\psi_{ttt}^0,&x\in\mb{R}^n,\ t>0,\\
		\psi^{I,1}(0,x)=\psi_t^{I,1}(0,x)=0,&x\in\mb{R}^n,
	\end{cases}
\end{align}
and
\begin{align}\label{Kuznetsov_I_2}
	\begin{cases}
		\psi_{tt}^{I,2}-\Delta\psi^{I,2}-\delta\Delta\psi_t^{I,2}=\Delta\psi_t^{I,1}-\psi_{ttt}^{I,1},&x\in\mb{R}^n,\ t>0,\\
		\psi^{I,2}(0,x)=\psi_t^{I,2}(0,x)=0,&x\in\mb{R}^n.
	\end{cases}
\end{align}
The deductions of the inhomogeneous strongly damped wave equations \eqref{Kuznetsov_I_1} and \eqref{Kuznetsov_I_2} will be shown precisely in Section \ref{Section-Singular-Limit} by using the WKB expansions and the multi-scale analysis.

\begin{theorem}\label{Theorem-Singular-Limit}Let $0<\tau\ll 1$. Suppose that the initial data $\psi_0\in H^{9+M_{p,q,n}}$ and $\psi_1\in H^{7+M_{p,q,n}}$ with $1\leqslant p\leqslant 2\leqslant q\leqslant +\infty$, $M_{p,q,n}$ defined in \eqref{Mpq} such that $\Delta^2(\psi_0+\delta\psi_1)=0$ as well as $\Delta\psi_2=0$. Then, the solution $\psi$ to the viscous MGT equation \eqref{Eq_MGT} satisfies the following $L^p-L^q$ refined estimates uniformly in time:
\begin{align*}
&\sup\limits_{t\in[0,+\infty)}\left\|\psi(t,\cdot)-\big(\psi^0(t,\cdot)+\tau\psi^{I,1}(t,\cdot)+\tau^2\psi^{I,2}(t,\cdot)\big)-\tau^2\psi^{\mathrm{Lay}}(t,\cdot)\right\|_{L^q}\\
&\leqslant C\tau^3\|(\psi_0,\psi_1)\|_{H^{9+M_{p,q,n}}_p\times H^{7+M_{p,q,n}}_p},
\end{align*}
 where the constant $C$ is independent of $\tau$, and the layer function $\psi^{\mathrm{Lay}}=\psi^{\mathrm{Lay}}(t,x)$ is denoted by
\begin{align}\label{fun-lay}
\psi^{\mathrm{Lay}}(t,x):=\left( \frac{t}{\tau}-1+\mathrm{e}^{-\frac{t}{\tau}}\right)\big(\psi_2(x)-\Delta\psi_0(x)-\delta\Delta\psi_1(x)\big).
\end{align}
\end{theorem}
\begin{remark}
Under the assumptions of Theorem \ref{Theorem-Singular-Limit}, we may have the global (in time) asymptotic expansion of the solution to the viscous MGT equation \eqref{Eq_MGT} as follows:
\begin{align*}
\psi(t,x)=\psi^0(t,x)+\tau\psi^{I,1}(t,x)+\tau^2\left(\psi^{I,2}(t,x)+\psi^{\mathrm{Lay}}(t,x)\right)+O(\tau^3)
\end{align*}
as $0<\tau\ll 1$, in the $L^{\infty}([0,+\infty),L^q)$ framework. Especially with $q=+\infty$, under suitable assumption for the data, we have shown
\begin{align*}
\psi\to\psi^0+\tau\psi^{I,1}+\tau^2\left(\psi^{I,2}+\psi^{\mathrm{Lay}}\right)\ \ \mbox{in}\ \ L^{\infty}([0,+\infty)\times\mb{R}^n)\ \ \mbox{as}\ \ \tau\downarrow0 
\end{align*}
with the rate of convergence $\tau^3$ for any $n\geqslant 1$.
 The layer function $\psi^{\mathrm{Lay}}$ with the fast change factor is one of the novelties in our manuscript. It improves the previous results \cite{B-C-Lasiecka=2020,Chen-Ikehta=2021,Chen-Takeda=2023} on the singular limit $\psi\to\psi^0$ with the linear rate $\tau$ by subtracting the additional error terms.
\end{remark}
\begin{remark}
	Concerning the small thermal relaxation profile of $\psi$ in Theorem \ref{Theorem-Singular-Limit}, we notice the fast change factor $\tau^2\psi^{\mathrm{Lay}}$ as $\tau\downarrow0$ whose dominant part can be expressed by $\Gamma(\tau;t_0)$ in Figure \ref{imggg} with a large fixed time $t=t_0\gg1$. The speed of the rapid change as $\tau\downarrow0$ depends on $t_0$, but in the exponential type in any case. Precisely, only when we consider $t\geqslant t_0\gg1$, the rapid change occurs (see the right figure).
	\begin{figure}[h]
		\centering
		\begin{tikzpicture}
			\draw[->] (-0.2,0) -- (5.8,0) node[below] {$\tau$};
			\draw[->] (0,-0.2) -- (0,2.8) node[left] {$\Gamma$};
			\node[left] at (0,-0.2) {{$0$}};
			\node[right] at (0.3,1) {$\leftarrow\Gamma(\tau;t_0):=\tau t_0-\tau^2+\tau^2\mathrm{e}^{-\frac{t_0}{\tau}}$};
			\draw[color=black] plot[smooth, tension=.7] coordinates {(0,0) (1,1.8) (5.4,2.1)};
			
			\draw[->] (7.8,0) -- (13.8,0) node[below] {$\tau$};
			\draw[->] (8,-0.2) -- (8,2.8) node[left] {$\Gamma$};
			\node[left] at (8,-0.2) {{$0$}};
			\node[below] at (11.3,2.8) {$\downarrow\Gamma(\tau;t_0)$};
			\draw[color=black] plot[smooth, tension=.7] coordinates {(8,0) (9,1.8) (13.4,2.1)};
			\draw[dashed, color=black] plot[smooth, tension=.7] coordinates {(8,0) (10,1) (13.4,1.2)};
			\node[below] at (11.3,0.8) {$\uparrow\Gamma(\tau;t_1)$\mbox{ with }$t_1\ll t_0$};
		\end{tikzpicture}
		\caption{A fast change as $\tau\downarrow0$ for a large fixed time $t=t_0\gg1$}
		\label{imggg}
	\end{figure}
\end{remark}

\begin{remark}
It is worth mentioning that we do not need to assume the compatible initial condition $\psi_2=\Delta(\psi_0+\delta\psi_1)$ since $\Delta\psi_2=\Delta^2(\psi_0+\delta\psi_1)$ implies
\begin{align*}
\psi_2=\Delta(\psi_0+\delta\psi_1)+\Upsilon\ \ \mbox{with}\ \ \Delta\Upsilon=0,
\end{align*}
where $\Upsilon=\Upsilon(x)$ is not necessary to be zero, e.g. $\Upsilon(x)=c_0+c_1x$.
\end{remark}

\begin{remark}
If $\psi_2=\Delta\psi_0+\delta\Delta\psi_1$, we find $\psi^{\mathrm{Lay}}\equiv0$. In other words, the fast change factor in the singular layer disappears due to the compatibility of initial conditions between \eqref{Eq_MGT} and \eqref{Eq_VEDW}, precisely, $\psi^0_{tt}(0,x)=\psi_{tt}(0,x)$. For this reason, it seems that the following identity:
\begin{align}\label{CON-01}
\psi_2(x)=\Delta\psi_0(x)+\delta\Delta\psi_1(x)\ \ \mbox{for}\ \ x\in\mb{R}^n,
\end{align}
is the threshold condition for the appearance of the initial layer for the viscous MGT equation \eqref{Eq_MGT}. We have partly verified this conjecture in Theorem \ref{Theorem-Singular-Limit}.
\end{remark}

\begin{remark}
	The appearance for the initial layer comes from the replacement of the Cattaneo law of heat conduction (hyperbolic-like version) by the Fourier law of heat conduction (parabolic-like version) for acoustic waves in thermoviscous flows.  This phenomenon also occurs between the classical damped wave equation and the heat equation (cf. \cite{Lions=1973,Ikehata=2003,Hashimoto-Yamazaki=2007,Ghisi-Gobbino=2012,Ikehata-Sobajima=2023} and references therein).
\end{remark}

\begin{remark}
In Theorem \ref{Theorem-Singular-Limit}, we rigorously justified the second order small thermal relaxation profiles to the viscous MGT equation \eqref{Eq_MGT} in $L^{\infty}([0,+\infty),L^q)$, namely, the correction for $j=2$ in the formal expansion proposed in Proposition \ref{Prop-WKB-formal}. We conjecture that the higher order profiles, that is, $j\geqslant 3$ in Proposition \ref{Prop-WKB-formal} also hold, which may be demonstrated by following the approaches in Section \ref{Section-Singular-Limit}.
\end{remark}

\section{Large time decay behavior}\label{Section-Lp-Lq}\setcounter{equation}{0}
\hspace{5mm}Unlike the previous researches \cite{Pellicer-Said-Houari=2019} via energy methods in the Fourier space and \cite{Chen-Ikehta=2021,Chen-Takeda=2023}  via explicit representations of solutions to the third order (in time) differential equation, we will study some large time qualitative properties of solutions to the viscous MGT equation \eqref{Eq_MGT} in the $L^q$ framework via diagonalization procedure, which gives the complete proof of Theorem \ref{Thm-Lp-Lq}. We may derive sharp behavior of solutions without using explicit expressions (see Remark \ref{Rem-Diag-Large-Class} for its motivation).

\subsection{The first order (in time) coupled system}
\hspace{5mm}To begin with our analysis, motivated by the wave-like structure of the viscous MGT equation \eqref{Eq_MGT}, namely,
\begin{align}\label{Motivation}
\underbrace{\psi_{tt}-\Delta\psi}_{\mbox{wave equation I}}+\partial_t\big(\underbrace{\tau\psi_{tt}-(\delta+\tau)\Delta\psi}_{\mbox{wave equation II}}\big)=0,
\end{align}
 let us introduce a suitable energy term $\Psi$ in \eqref{Energy-Term-Psi}. The propagation speed of the wave equation II, that is $c_0=\sqrt{\frac{\delta+\tau}{\tau}}$, decides the dominant speed of acoustic waves propagation for the viscous MGT equation \eqref{Eq_MGT}, which has been noted in Remark \ref{Rem-Prop-sing} and \cite[Remark 3.4]{Chen-Ikehta=2021}. The choice of a suitable energy term plays the most important role in our diagonalization procedure. It solves the first order (in time) coupled system
\begin{align}\label{Eq-Psi}
\begin{cases}
\Psi_t+(A_0+A_1|D|)\Psi=0,&x\in\mb{R}^n,\ t>0,\\
\Psi(0,x)=\Psi_0(x),&x\in\mb{R}^n,
\end{cases}
\end{align}
where the coefficient matrices are given by
\begin{align*}
	A_0:=\frac{1}{\tau}\left(\begin{array}{ccc}
		0 & 0 & -1\\
		0 & 0 & -1\\
		0 & 0 & 1
	\end{array}\right)\ \ \mbox{and}\ \ A_1:=\frac{1}{2\tau}\left(\begin{array}{ccc}
		\delta & -(2\tau+\delta) & 0 \\
		2\tau+\delta & -\delta & 0\\ 
		-\delta & \delta & 0
	\end{array}\right).
\end{align*}
\begin{remark}
The classical and general results for hyperbolic(-parabolic) coupled systems do not work well in our model \eqref{Eq-Psi}. For instance, the non-symmetric matrices $A_0$ and $A_1$ violate the condition in \cite[Lemma 2.2]{Umeda-Kawashima-Shizuta}. Later, strongly motivated by diagonalization procedure developed in \cite{Reissig-Wang=2005,Yang-Wang=2006,Jachmann=2008}, we will decouple the system \eqref{Eq-Psi} by the frequency analysis to study sharp $L^p-L^q$ decay estimates and large time asymptotic profile of solutions. This philosophy can be widely applied in the MGT-type coupled systems in the future.
\end{remark}

Let us apply the partial Fourier transform with respect to spatial variables for the Cauchy problem \eqref{Eq-Psi}, we may deduce
\begin{align}\label{Psi-Fourier}
\begin{cases}
\widehat{\Psi}_t+(A_0+A_1|\xi|)\widehat{\Psi}=0,&\xi\in\mb{R}^n,\ t>0,\\
\widehat{\Psi}(0,\xi)=\widehat{\Psi}_0(\xi),&\xi\in\mb{R}^n.
\end{cases}
\end{align}
Before applying the WKB analysis, we take the next zones in the Fourier space:
\begin{align*}
	\ml{Z}_{\intt}(\varepsilon_0):=\{|\xi|\leqslant\varepsilon_0\ll1\},\ \ 
	\ml{Z}_{\bdd}(\varepsilon_0,N_0):=\{\varepsilon_0\leqslant |\xi|\leqslant N_0\},\ \  
	\ml{Z}_{\extt}(N_0):=\{ |\xi|\geqslant N_0\gg1\}.
\end{align*}
The cut-off functions $\chi_{\intt}(\xi),\chi_{\bdd}(\xi),\chi_{\extt}(\xi)\in \mathcal{C}^{\infty}$ having supports in their corresponding zones $\ml{Z}_{\intt}(\varepsilon_0)$, $\ml{Z}_{\bdd}(\varepsilon_0/2,2N_0)$ and $\ml{Z}_{\extt}(N_0)$, respectively, satisfying
\begin{align*}
\chi_{\bdd}(\xi)=1-\chi_{\intt}(\xi)-\chi_{\extt}(\xi).
\end{align*}
Our consideration for deriving asymptotic behavior of $\Psi$ in the Fourier space is divided into the next three subclasses in regard to the magnitude of frequencies:
\begin{enumerate}[(1)]
\item asymptotic representation of $\widehat{\Psi}$ when $\xi\in\ml{Z}_{\intt}(\varepsilon_0)$;
\item asymptotic representation of $\widehat{\Psi}$ when $\xi\in\ml{Z}_{\extt}(N_0)$;
\item exponential stability of $\widehat{\Psi}$ when $\xi\in\ml{Z}_{\bdd}(\varepsilon_0,N_0)$.
\end{enumerate}

\subsection{Subclass I: Diagonalization procedure for small frequencies}\label{Sub-Section-Diag-Small}
\hspace{5mm}Due to the coefficient $A_0$ exerting dominant influence comparing with another coefficient $A_1|\xi|$ for $\xi\in\ml{Z}_{\intt}(\varepsilon_0)$, let us begin with our diagonalization procedure for the matrix $A_0$. We introduce the first quantity $\widehat{\Psi}^{(1,s)}:=T_{1,s}^{-1}\widehat{\Psi}$ with the matrix $T_{1,s}$ defined in \eqref{T1,N2,T3}, which may diagonalize $A_0$. As a consequence, multiplying $T_{1,s}^{-1}$ on the equation of the Cauchy problem \eqref{Psi-Fourier}, this quantity solves
\begin{align}\label{Eq-Psi-1}
\widehat{\Psi}^{(1,s)}_t+\Lambda_1^{(s)}\widehat{\Psi}^{(1,s)}+R_1^{(s)}|\xi|\widehat{\Psi}^{(1,s)}=0,
\end{align}
whose coefficient matrices are addressed by
\begin{align*}
\Lambda_1^{(s)}=\frac{1}{\tau}\diag(0,0,1)\ \ \mbox{and}\ \ 
R_1^{(s)}=\frac{1}{2\tau}\left(\begin{array}{ccc}
	0 & -2\tau & 2\tau\\
	2\tau & 0 & -2\tau\\ 
	-\delta & \delta & 0
\end{array}\right).
\end{align*}

With the purpose of retaining the derived diagonal part $\Lambda_1^{(s)}$ in \eqref{Eq-Psi-1}, we may introduce the second quantity $\widehat{\Psi}^{(2,s)}:=T_{2,s}^{-1}\widehat{\Psi}^{(1,s)}$, where we defined $T_{2,s}:=I_{3}+N_{2,s}|\xi|$ with the auxiliary matrix $N_{2,s}$ showed in \eqref{T1,N2,T3}. Let us recall $T^{-1}=I_{3}-NT^{-1}$ for $T=I_{3}+N$.  According to the facts that
\begin{align}\label{Exp-01}
T_{2,s}^{-1}\Lambda_1^{(s)}T_{2,s}&=\Lambda_1^{(s)}-N_{2,s}T_{2,s}^{-1}\Lambda_1^{(s)}T_{2,s}|\xi|+\Lambda_1^{(s)}N_{2,s}|\xi|\notag\\
&=\Lambda_1^{(s)}-[N_{2,s},\Lambda_1^{(s)}]|\xi|-N_{2,s}\Lambda_1^{(s)}N_{2,s}|\xi|^2+N_{2,s}^2T_{2,s}^{-1}\Lambda_1^{(s)}T_{2,s}|\xi|^2,
\end{align}
and similarly,
\begin{align}\label{Exp-02}
T_{2,s}^{-1}R_1^{(s)}T_{2,s}|\xi|=R_1^{(s)}|\xi|-N_{2,s}T_{2,s}^{-1}R_1^{(s)}T_{2,s}|\xi|^2+R_1^{(s)}N_{2,s}|\xi|^2,
\end{align}
by directly multiplying \eqref{Eq-Psi-1} by $T_{2,s}^{-1}$, we are able to find
\begin{align}\label{Eq-Psi-2}
\widehat{\Psi}^{(2,s)}_t+\Lambda_1^{(s)}\widehat{\Psi}^{(2,s)}+\left(R_1^{(s)}-[N_{2,s},\Lambda_1^{(s)}]\right)|\xi|\widehat{\Psi}^{(2,s)}+R_2^{(s)}|\xi|^2\widehat{\Psi}^{(2,s)}=0,
\end{align}
where $[A,B]:=AB-BA$ denotes the commutator of matrices and the remainder is expressed by
\begin{align*}
R_2^{(s)}=-N_{2,s}\Lambda_1^{(s)}N_{2,s}+N_{2,s}^2T_{2,s}^{-1}\Lambda_1^{(s)}T_{2,s}-N_{2,s}T_{2,s}^{-1}R_1^{(s)}T_{2,s}+R_1^{(s)}N_{2,s}.
\end{align*}
Thanks to the construction of $N_{2,s}$ in \eqref{T1,N2,T3}, it comes
\begin{align*}
A_2^{(s)}=A_{1}^{(s)}-[N_{2,s},\Lambda_1^{(s)}]=\left(\begin{array}{cc;{2pt/2pt}c}
	0 & -1 & 0\\
	1 & 0 & 0\\ \hdashline[2pt/2pt] 
	0 & 0 & 0
\end{array}\right).
\end{align*}
In other words, we just can get the block diagonal matrix $A_2^{(s)}$ so far. To achieve our aim of deriving the next diagonal part with the size $|\xi|$, let us take $\widehat{\Psi}^{(3,s)}:=T_{3,s}^{-1}\widehat{\Psi}^{(2,s)}$ with the matrix $T_{3,s}$ defined in \eqref{T1,N2,T3}, which may diagonalize $A_2^{(s)}$ and retain the first diagonal matrix $\Lambda_1^{(s)}$ simultaneously, because one may regard $\Lambda_1^{(s)}$ as the block diagonal matrix also.  Multiplying \eqref{Eq-Psi-2} by $T_{3,s}^{-1}$, we immediately arrive at
\begin{align*}
\widehat{\Psi}^{(3,s)}_t+\left(\Lambda_{1}^{(s)}+\Lambda_2^{(s)}|\xi|\right)\widehat{\Psi}^{(3,s)}+R_3^{(s)}|\xi|^2\widehat{\Psi}^{(3,s)}=0
\end{align*}
with the coefficient matrices $\Lambda_2^{(s)}=\diag(i,-i,0)$ and $R_3^{(s)}=T_{3,s}^{-1}R_2^{(s)}T_{3,s}$.

Before carrying out the next step, we extract the lowest order term $R_{3,1}^{(s)}$ to be the dominant term  of the remainder $R_3^{(s)}$ as follows:
\begin{align*}
R_3^{(s)}&=T_{3,s}^{-1}\left(N_{2,s}[N_{2,s},\Lambda_1^{(s)}]+[R_1^{(s)},N_{2,s}]\right)T_{3,s}\\
&\quad+T_{3,s}^{-1}\left(N_{2,s}^2\Lambda_1^{(s)}N_{2,s}-N_{2,s}^3T_{2,s}^{-1}\Lambda_1^{(s)}T_{2,s}+N_{2,s}^{2}T_{2,s}^{-1}R_1^{(s)}T_{2,s}-N_{2,s}R_1^{(s)}N_{2,s}\right)T_{3,s}|\xi|\\
&=:R_{3,1}^{(s)}+R_{3,2}^{(s)}|\xi|,
\end{align*}
particularly,
\begin{align*}
R_{3,1}^{(s)}=\frac{1}{2}\left(\begin{array}{cc;{2pt/2pt}c}
	\delta & \delta i & (1-i)\tau\\
	-\delta i & \delta & (1+i)\tau\\ \hdashline[2pt/2pt] 
	(1+i)\delta & (1-i)\delta & -2\delta
\end{array}\right).
\end{align*}
To determine the third diagonal matrix, we divide our discussion into two substeps.\\
\textbf{\underline{Substep I}:} Because identical value of two diagonal elements of $\Lambda_1^{(s)}$, the first left block matrix of $R_{3,1}^{(s)}$ cannot be fully diagonalized by $\Lambda_1^{(s)}$ in the first substep. Let us define $\widehat{\Psi}^{(3\frac{1}{2},s)}:=T^{-1}_{3\frac{1}{2},s}\widehat{\Psi}^{(3,s)}$, where we set $T_{3\frac{1}{2},s}:=I_{3}+N_{3\frac{1}{2},s}|\xi|^2$ with the auxiliary matrix $N_{3\frac{1}{2},s}$ defined in \eqref{N3-1/2,N4}.
 With the same manners of expressions in \eqref{Exp-01} and \eqref{Exp-02}, we are able to derive
\begin{align*}
	\widehat{\Psi}^{(3\frac{1}{2},s)}_t+\left(\Lambda_{1}^{(s)}+\Lambda_2^{(s)}|\xi|+A_3^{(s)}|\xi|^2\right)\widehat{\Psi}^{(3\frac{1}{2},s)}+R_{3\frac{1}{2}}^{(s)}|\xi|^3\widehat{\Psi}^{(3\frac{1}{2},s)}=0,
\end{align*}
in which the coefficient matrices are given by
\begin{align*}
A_3^{(s)}&=R_{3,1}^{(s)}-[N_{3\frac{1}{2},s},\Lambda_1^{(s)}]=\frac{1}{2}\left(\begin{array}{cc;{2pt/2pt}c}
	\delta & \delta i & 0\\
	-\delta i & \delta & 0\\ \hdashline[2pt/2pt] 
	0 & 0 & -2\delta
\end{array}\right),\\
R_{3\frac{1}{2}}^{(s)}&=T_{3\frac{1}{2},s}^{-1}R_{3,2}^{(s)}T_{3\frac{1}{2},s}-N_{3\frac{1}{2},s}\Lambda_1^{(s)}N_{3\frac{1}{2},s}|\xi|+N_{3\frac{1}{2},s}^2T_{3\frac{1}{2},s}^{-1}\Lambda_1^{(s)}T_{3\frac{1}{2},s}|\xi|-N_{3\frac{1}{2},s}T_{3\frac{1}{2},s}^{-1}\Lambda_2^{(s)}T_{3\frac{1}{2},s}\\
&\quad+\Lambda_2^{(s)}N_{3\frac{1}{2},s}-N_{3\frac{1}{2},s}T_{3\frac{1}{2},s}^{-1}R_{3,1}^{(s)}T_{3\frac{1}{2},s}|\xi|+R_{3,1}^{(s)}N_{3\frac{1}{2},s}|\xi|.
\end{align*}
\textbf{\underline{Substep II}:} Let us take additionally $\widehat{\Psi}^{(4,s)}:=T_{4,s}^{-1}\widehat{\Psi}^{(3\frac{1}{2},s)}$ carrying $T_{4,s}:=I_{3}+N_{4,s}|\xi|$ with the auxiliary matrix $N_{4,s}$ defined in \eqref{N3-1/2,N4}. This construction is motivated by getting a diagonal part of the first left block matrix of $A_3^{(s)}$. Thanks to $T_{4,s}^{-1}\Lambda_1^{(s)}T_{4,s}=\Lambda_1^{(s)}$, we claim
\begin{align*}
\widehat{\Psi}_t^{(4,s)}+\left(\Lambda_1^{(s)}+\Lambda_{2}^{(s)}|\xi|+\Lambda_3^{(s)}|\xi|^2\right)\widehat{\Psi}^{(4,s)}+R_4^{(s)}|\xi|^3\widehat{\Psi}^{(4,s)}=0,
\end{align*}
where two coefficient matrices are shown via
\begin{align*}
\Lambda_3^{(s)}&=A_{3}^{(s)}-[N_{4,s},\Lambda_2^{(s)}]=\frac{1}{2}\diag(\delta,\delta,-2\delta),\\
R_4^{(s)}&=T_{4,s}^{-1}R_{3\frac{1}{2}}^{(s)}T_{4,s}-N_{4,s}\Lambda_2^{(s)}N_{4,s}+N_{4,s}^2T_{4,s}^{-1}\Lambda_2^{(s)}T_{4,s}-N_{4,s}T_{4,s}^{-1}A_3^{(s)}T_{4,s}+A_3^{(s)}N_{4,s}.
\end{align*}

From the previous diagonalization procedure, we have obtained pairwise distinct characteristic roots from $\Lambda_1^{(s)}+\Lambda_{2}^{(s)}|\xi|+\Lambda_3^{(s)}|\xi|^2$ whose real parts are strictly negative. Consequently, following the general philosophy proposed in \cite[Chapter 2]{Jachmann=2008}, we may derive asymptotic representation of the  energy term for small frequencies.
\begin{prop}\label{Prop-small-frenquencies}
	The characteristic roots $\lambda_j(|\xi|)$ of the coefficient matrix $A_0+A_1|\xi|$ in the Cauchy problem \eqref{Psi-Fourier} behave for $\xi\in\ml{Z}_{\intt}(\varepsilon_0)$ as
	\begin{align*}
		\lambda_{1,2}(|\xi|)=\mp i|\xi|-\frac{\delta}{2}|\xi|^2+O(|\xi|^3), \ \ 	\lambda_{3}(|\xi|)=-\frac{1}{\tau}+\delta|\xi|^2+O(|\xi|^3).
	\end{align*}
	Notice that these expansions coincide with those in the context of \cite{Pellicer-Said-Houari=2019,Chen-Ikehta=2021,Chen-Takeda=2023}.
	The solution to the Cauchy problem \eqref{Psi-Fourier} has the asymptotic representation for $\xi\in\ml{Z}_{\intt}(\varepsilon_0)$ as follows: 
	\begin{align*}   
		\chi_{\intt}(\xi)\widehat{\Psi}(t,\xi)=\chi_{\intt}(\xi)T_{\intt}\diag\big(\mathrm{e}^{\lambda_j(|\xi|)t}\big)_{j=1}^3 T_{\intt}^{-1}\widehat{\Psi}_0(\xi)
	\end{align*}
	with $T_{\intt}:=T_{1,s}(I_{3}+N_{2,s}|\xi|)T_{3,s}(I_{3}+N_{3\frac{1}{2},s}|\xi|^2)(I_{3}+N_{4,s}|\xi|)$, where these auxiliary matrices are defined in \eqref{T1,N2,T3} as well as \eqref{N3-1/2,N4}.
\end{prop}

\subsection{Subclass II: Diagonalization procedure for large frequencies}
\hspace{5mm}Let us now turn to the subclass for $\xi\in\ml{Z}_{\extt}(N_0)$, in which the coefficient $A_1|\xi|$ plays the dominant role in comparison with another matrix $A_0$. For this reason, to diagonalize $A_1$ directly, we introduce the new quantity $\widehat{\Psi}^{(1,l)}:=T_{1,l}^{-1}\widehat{\Psi}$ with the matrix
\begin{align}\label{T1l}
	T_{1,l}:=\frac{1}{\delta}\left(\begin{array}{ccc}
		0 & -(\delta+\tau)-i\sqrt{\tau(\delta+\tau)} & -(\delta+\tau)+i\sqrt{\tau(\delta+\tau)}\\
		0 & -(\delta+\tau)+i\sqrt{\tau(\delta+\tau)} & -(\delta+\tau)-i\sqrt{\tau(\delta+\tau)}\\ 
		\delta & \delta & \delta
	\end{array}\right).
\end{align}
Hence, it fulfills
\begin{align*}
\widehat{\Psi}^{(1,l)}_t+\Lambda_1^{(l)}|\xi|\widehat{\Psi}^{(1,l)}+R_1^{(l)}\widehat{\Psi}^{(1,l)}=0
\end{align*}
with the coefficient matrices
\begin{align*}
\Lambda_1^{(l)}=\frac{1}{\sqrt{\tau}}\diag\left(0,i\sqrt{\delta+\tau},-i\sqrt{\delta+\tau}\,\right)\ \ \mbox{and}\ \ R_1^{(l)}=\frac{1}{2\tau(\delta+\tau)}\left(\begin{array}{ccc}
	2\tau & 2\tau & 2\tau\\
	\delta & \delta & \delta\\ 
	\delta & \delta & \delta
\end{array}\right).
\end{align*}

Analogously to the subclass for small frequencies, let us consider $\widehat{\Psi}^{(2,l)}:=T_{2,l}^{-1}\widehat{\Psi}^{(1,l)}$ with $T_{2,l}:=I_{3}+N_{2,l}|\xi|^{-1}$ carrying the auxiliary matrix
\begin{align}\label{N2l}
N_{2,l}:=\frac{i}{4(\delta+\tau)\sqrt{\tau(\delta+\tau)}}\left(\begin{array}{ccc}
	0 & -4\tau  & 4\tau \\
	2\delta & 0 & \delta \\ 
	-2\delta & -\delta & 0
\end{array}\right).
\end{align}
By the same techniques as those in \eqref{Exp-01} as well as \eqref{Exp-02}, we notice that
\begin{align*}
\widehat{\Psi}^{(2,l)}_t+\left(\Lambda_1^{(l)}|\xi|+\Lambda_2^{(l)}\right)\widehat{\Psi}^{(2,l)}+R_2^{(l)}|\xi|^{-1}\widehat{\Psi}^{(2,l)}=0,
\end{align*}
in which the crucial diagonal part is
\begin{align*}
\Lambda_2^{(l)}=R_1^{(l)}-[N_{2,l},\Lambda_1^{(l)}]=\frac{1}{2\tau(\delta+\tau)}\diag(2\tau,\delta,\delta),
\end{align*}
moreover, the remainder is
\begin{align*}
R_2^{(l)}=-N_{2,l}\Lambda_1^{(l)}N_{2,l}+N_{2,l}^2T_{2,l}^{-1}\Lambda_1^{(l)}T_{2,l}-N_{2,l}T_{2,l}^{-1}R_1^{(l)}T_{2,l}+R_1^{(l)}N_{2,l}.
\end{align*}

According to the previous diagonalization procedure, we have derived pairwise distinct characteristic roots from $\Lambda_1^{(l)}|\xi|+\Lambda_{2}^{(l)}$ whose real parts are strictly negative. Then, following the general philosophy proposed in \cite[Chapter 2]{Jachmann=2008} again, we may obtain asymptotic representation of the  energy term for large frequencies.
\begin{prop}\label{Prop-large-frenquencies}
	The characteristic roots $\lambda_j(|\xi|)$ of the coefficient matrix $A_0+A_1|\xi|$ in the Cauchy problem \eqref{Psi-Fourier} behave for $\xi\in\ml{Z}_{\extt}(N_0)$ as
	\begin{align*}
		\lambda_{1}(|\xi|)=-\frac{1}{\delta+\tau}+O(|\xi|^{-1}), \ \ \lambda_{2,3}(|\xi|)=\pm i\sqrt{\frac{\delta+\tau}{\tau}}|\xi|-\frac{\delta}{2\tau(\delta+\tau)}+O(|\xi|^{-1}).
	\end{align*}
	Notice that these expansions coincide with those in the context of \cite{Pellicer-Said-Houari=2019,Chen-Ikehta=2021,Chen-Takeda=2023}. The solution to the Cauchy problem \eqref{Psi-Fourier} has the asymptotic representation for $\xi\in\ml{Z}_{\extt}(N_0)$ as follows: 
	\begin{align*}   
		\chi_{\extt}(\xi)\widehat{\Psi}(t,\xi)=\chi_{\extt}(\xi)T_{\extt}\diag\big(\mathrm{e}^{\lambda_j(|\xi|)t}\big)_{j=1}^3 T_{\extt}^{-1}\widehat{\Psi}_0(\xi)
	\end{align*}
	with $T_{\extt}:=T_{1,l}(I_{3}+N_{2,l}|\xi|^{-1})$, where these auxiliary matrices are defined in \eqref{T1l} and \eqref{N2l}.
\end{prop}

\subsection{Subclass III: Exponential stability for bounded frequencies}
\hspace{5mm}Finally, we investigate an exponential decay result for the energy term as $\xi\in\ml{Z}_{\bdd}(\varepsilon_0,N_0)$ by employing a contradiction argument. Suppose that there is a pure imaginary characteristic root, that is $\exists j\in\{1,2,3\}$ such that $\lambda_j(|\xi|)=ia(|\xi|)$ with $a=a(|\xi|)\in\mb{R}\backslash\{0\}$, of the coefficient matrix $A_0+A_1|\xi|$ for $\xi\in\ml{Z}_{\bdd}(\varepsilon_0,N_0)$. The non-trivial $|\xi|$-dependent real number $a$ satisfies the following identity:
\begin{align*}
0&=\mathrm{det}(A_0+A_1|\xi|-iaI_{3})=\left|\begin{array}{ccc}
	\frac{\delta}{2\tau}|\xi|-ia & -\frac{2\tau+\delta}{2\tau}|\xi|  & -\frac{1}{\tau} \\[0.28em]
	\frac{2\tau+\delta}{2\tau}|\xi| & -\frac{\delta}{2\tau}|\xi|-ia & -\frac{1}{\tau} \\ [0.28em]
	-\frac{\delta}{2\tau}|\xi| & \frac{\delta}{2\tau}|\xi| & \frac{1}{\tau}-ia
\end{array}\right|\\
&=ia\left(a^2-\frac{\delta+\tau}{\tau}|\xi|^2\right)-\frac{1}{\tau}(a^2-|\xi|^2).
\end{align*}

Collecting the real part and imaginary part separately, we find that $a^2=\frac{\delta+\tau}{\tau}|\xi|^2$ and $a^2=|\xi|^2$, which cannot be achieved simultaneously due to $\delta>0$ in the viscous MGT equation. Namely, it gives a contradiction so that the characteristic roots of $A_0+A_1|\xi|$ cannot be pure imaginary numbers. From the compactness of the bounded zone $\ml{Z}_{\bdd}(\varepsilon_0,N_0)$ and the continuity of $\lambda_j(|\xi|)$ in regard to $|\xi|$, we conclude $\mathrm{Re}\,\lambda_j(|\xi|)<0$ for $\xi\in\ml{Z}_{\bdd}(\varepsilon_0,N_0)$ and the next exponentially stable result.
\begin{prop}\label{Prop-bounded-frenquencies}
	The solution to the Cauchy problem \eqref{Psi-Fourier} fulfills an exponential decay estimate
	\begin{align*}
		\chi_{\bdd}(\xi)|\widehat{\Psi}(t,\xi)|\lesssim \chi_{\bdd}(\xi)\mathrm{e}^{-ct}|\widehat{\Psi}_0(\xi)|
	\end{align*}
with a positive constant $c$.
\end{prop}

\subsection{Decay property and large time profile: Proof of Theorem \ref{Thm-Lp-Lq}}\label{Sub-Est}
\hspace{5mm}We firstly derive the $L^p-L^q$ estimates when $1\leqslant p\leqslant 2\leqslant q\leqslant+\infty$. 
 Let us apply Proposition \ref{Prop-small-frenquencies} associated with the Hausdorff-Young inequality and H\"older's inequality to obtain
\begin{align*}
\|\chi_{\intt}(D)|D|^s\Psi(t,\cdot)\|_{(L^q)^3}&\lesssim\|\chi_{\intt}(\xi)|\xi|^s\mathrm{e}^{-c|\xi|^2t}\widehat{\Psi}_0(\xi)\|_{(L^{q'})^3}\\
&\lesssim\begin{cases}
\|\chi_{\intt}(\xi)|\xi|^s\mathrm{e}^{-c|\xi|^2t}\|_{L^{\frac{p'q'}{p'-q'}}}\|\widehat{\Psi}_0\|_{(L^{p'})^3}&\mbox{when} \ \ p\neq q,\\
\|\chi_{\intt}(\xi)|\xi|^s\mathrm{e}^{-c|\xi|^2t}\|_{L^{\infty}}\|\widehat{\Psi}_0\|_{(L^{p'})^3}&\mbox{when} \ \ p= q,
\end{cases}
\end{align*}
with $\frac{1}{p}+\frac{1}{p'}=1=\frac{1}{q}+\frac{1}{q'}$. According to the inequalities that
\begin{align*}
\|\chi_{\intt}(\xi)|\xi|^s\mathrm{e}^{-c|\xi|^2t}\|_{L^{\frac{p'q'}{p'-q'}}}^{\frac{p'q'}{p'-q'}}
&\lesssim\int_{|\xi|\leqslant\varepsilon_0}|\xi|^{\frac{spq}{q-p}}\mathrm{e}^{-\frac{cpq}{q-p}|\xi|^2t}\mathrm{d}\xi\lesssim(1+t)^{-\frac{spq}{2(q-p)}-\frac{n}{2}}\ \ \mbox{when} \ \ p\neq q,\\
\|\chi_{\intt}(\xi)|\xi|^s\mathrm{e}^{-c|\xi|^2t}\|_{L^{\infty}}
&\lesssim (1+t)^{-\frac{s}{2}}\ \ \mbox{when}\ \ p=q,
\end{align*}
we immediately conclude
\begin{align*}
\|\chi_{\intt}(D)|D|^s\Psi(t,\cdot)\|_{(L^q)^3}\lesssim (1+t)^{-\frac{s}{2}-\frac{n}{2}(\frac{1}{p}-\frac{1}{q})}\|\Psi_0\|_{(L^{p})^3}.
\end{align*}
Analogously, we next use Proposition \ref{Prop-large-frenquencies} and find
\begin{align*}
\|\chi_{\extt}(D)|D|^s\Psi(t,\cdot)\|_{(L^q)^3}&\lesssim\|\chi_{\extt}(\xi)\langle\xi\rangle^s\mathrm{e}^{-ct}\widehat{\Psi}_0(\xi)\|_{(L^{q'})^3}\\
&\lesssim\begin{cases}
\mathrm{e}^{-ct}\|\chi_{\extt}(\xi)\langle\xi\rangle^{-M_{p,q,n}}\|_{L^{\frac{p'q'}{p'-q'}}}\|\langle\xi\rangle^{s+M_{p,q,n}}\widehat{\Psi}_0\|_{(L^{p'})^3}&\mbox{when}\ \ p\neq q,\\
\mathrm{e}^{-ct}\|\langle\xi\rangle^{s}\widehat{\Psi}_0\|_{(L^{p'})^3}&\mbox{when}\ \ p= q.
\end{cases}
\end{align*}
Hereafter, we denote the Japanese bracket by $\langle\xi\rangle^2=1+|\xi|^2$. Due to the chain that
\begin{align}\label{Large-Infin}
\|\chi_{\extt}(\xi)\langle\xi\rangle^{-M_{p,q,n}}\|_{L^{\frac{p'q'}{p'-q'}}}^{\frac{p'q'}{p'-q'}}\lesssim\int_{|\xi|\geqslant N_0}\langle\xi\rangle^{-M_{p,q,n}\frac{pq}{q-p}}\mathrm{d}\xi\lesssim\int_{N_0}^{+\infty}\langle r\rangle^{n-1-M_{p,q,n}\frac{pq}{q-p}}\mathrm{d}r<+\infty
\end{align}
when $M_{p,q,n}>n(\frac{1}{p}-\frac{1}{q})$ as $p\neq q$, one arrives at
\begin{align*}
\|\chi_{\extt}(D)|D|^s\Psi(t,\cdot)\|_{(L^q)^3}\lesssim \mathrm{e}^{-ct}\|\Psi_0\|_{(H^{s+M_{p,q,n}}_p)^3}.
\end{align*}
Finally,  because of an exponential decay for bounded frequencies, summarizing the last derived estimates, we may have
\begin{align*}
\|\Psi(t,\cdot)\|_{(\dot{H}^s_q)^3}&\lesssim(1+t)^{-\frac{s}{2}-\frac{n}{2}(\frac{1}{p}-\frac{1}{q})}\|\Psi_0\|_{(L^{p})^3}+\mathrm{e}^{-ct}\|\Psi_0\|_{(L^p)^3}+ \mathrm{e}^{-ct}\|\Psi_0\|_{(H^{s+M_{p,q,n}}_p)^3}\\
&\lesssim  (1+t)^{-\frac{s}{2}-\frac{n}{2}(\frac{1}{p}-\frac{1}{q})}\|\Psi_0\|_{(H^{s+M_{p,q,n}}_p)^3}
\end{align*}
for any $1\leqslant p\leqslant 2\leqslant q\leqslant+\infty$ and $M_{p,q,n}$ already defined in \eqref{Mpq}. This completes the first part of the theorem. As a byproduct, according to Propositions \ref{Prop-small-frenquencies}-\ref{Prop-bounded-frenquencies}, the well-posedness result can be demonstrated straightforwardly.

Let us focus on large time asymptotic profile with some refined estimates in the $L^q$ framework. To begin with our deduction, we introduce the $|\xi|$-dependent numbers $\mu_j(|\xi|)$ such that
\begin{align*}
\mu_{1,2}(|\xi|):=\mp i|\xi|-\frac{\delta}{2}|\xi|^2,\ \ \mu_3(|\xi|):=-\frac{1}{\tau}+\delta|\xi|^2.
\end{align*}
By denoting $\widetilde{T}_{\intt}:=T_{\intt}T_{4,s}^{-1}$, we recall the reference system \eqref{Ref-System} carrying its initial data in the Fourier space $\widehat{\Phi}_0(\xi)=\widetilde{T}_{\intt}^{-1}\widehat{\Psi}_0(\xi)$. For this reason, its solution is uniquely given by
\begin{align*}
\widehat{\Phi}(t,\xi)=\diag\left(\mathrm{e}^{\mu_j(|\xi|)t}\right)_{j=1}^3\widetilde{T}_{\intt}^{-1}\widehat{\Psi}_0(\xi).
\end{align*}
Making use of
\begin{align*}
T_{4,s}=I_{3}+N_{4,s}|\xi|\ \ \mbox{and}\ \ T_{4,s}^{-1}=I_{3}-N_{4,s}T_{4,s}^{-1}|\xi|,
\end{align*}
we may decompose the solution $\widehat{\Psi}$ for small frequencies into three parts as follows:
\begin{align*}
\chi_{\intt}(\xi)\widehat{\Psi}(t,\xi)&=\chi_{\intt}(\xi)\widetilde{T}_{\intt}T_{4,s}\diag\left(\mathrm{e}^{\lambda_j(|\xi|)t}\right)_{j=1}^3T_{4,s}^{-1}\widetilde{T}_{\intt}^{-1}\widehat{\Psi}_0(\xi)\\
&=\chi_{\intt}(\xi)\widetilde{T}_{\intt}\diag\left(\mathrm{e}^{\lambda_j(|\xi|)t}\right)_{j=1}^3\widetilde{T}_{\intt}^{-1}\widehat{\Psi}_0(\xi)\\
&\quad+\chi_{\intt}(\xi)|\xi|\widetilde{T}_{\intt}N_{4,s} \diag\left(\mathrm{e}^{\lambda_j(|\xi|)t}\right)_{j=1}^3T_{4,s}^{-1}\widetilde{T}_{\intt}^{-1}\widehat{\Psi}_0(\xi)\\
&\quad-\chi_{\intt}(\xi)|\xi|\widetilde{T}_{\intt} \diag\left(\mathrm{e}^{\lambda_j(|\xi|)t}\right)_{j=1}^3N_{4,s}T_{4,s}^{-1}\widetilde{T}_{\intt}^{-1}\widehat{\Psi}_0(\xi)\\
&=:S_{\intt,1}(t,\xi)+S_{\intt,2}(t,\xi)+S_{\intt,3}(t,\xi).
\end{align*}
One finds the refined estimate
\begin{align*}
|S_{\intt,1}(t,\xi)-\chi_{\intt}(\xi)\widetilde{T}_{\intt}\widehat{\Phi}(t,\xi)|&\lesssim\chi_{\intt}(\xi)\sup\limits_{j=1,2,3}\left|\mathrm{e}^{\lambda_j(|\xi|)t}-\mathrm{e}^{\mu_j(|\xi|)t}\right||\widehat{\Psi}_0(\xi)|\\
&\lesssim\chi_{\intt}(\xi)|\xi|(|\xi|^2t)\mathrm{e}^{-c|\xi|^2t}\int_0^1\mathrm{e}^{O(|\xi|^3)t\eta}\mathrm{d}\eta|\widehat{\Psi}_0(\xi)|\\
&\lesssim\chi_{\intt}(\xi)|\xi|\mathrm{e}^{-c|\xi|^2t}|\widehat{\Psi}_0(\xi)|.
\end{align*}
An application of H\"older's inequality implies
\begin{align*}
\left\|\chi_{\intt}(\xi)|\xi|^s\big(\widehat{\Psi}(t,\xi)-\widetilde{T}_{\intt}\widehat{\Phi}(t,\xi)\big)\right\|_{(L^{q'})^3}&\lesssim\left\|\chi_{\intt}(\xi)|\xi|^{s+1}\mathrm{e}^{-c|\xi|^2t}\right\|_{L^{\frac{p'q'}{p'-q'}}}\|\widehat{\Psi}_0\|_{(L^{p'})^3}\\
&\lesssim (1+t)^{-\frac{1}{2}-\frac{s}{2}-\frac{n}{2}(\frac{1}{p}-\frac{1}{q})}\|\Psi_0\|_{(L^p)^3}.
\end{align*}
Concerning the other parts, we repeat the investigation of \eqref{Large-Infin} to claim
\begin{align*}
\left\|\big(1-\chi_{\intt}(\xi)\big)|\xi|^s\widetilde{T}_{\intt}\widehat{\Phi}(t,\xi)\right\|_{(L^{q'})^3}&\lesssim\left\|\big(1-\chi_{\intt}(\xi)\big)\langle\xi\rangle^s\mathrm{e}^{-c|\xi|^2t}\widehat{\Psi}_0(\xi)\right\|_{(L^{q'})^3}\\
&\lesssim\mathrm{e}^{-ct}\|\Psi_0\|_{(H_p^{s+M_{p,q,n}})^3}.
\end{align*}
Summarizing the above obtained estimates and  using the Hausdorff-Young inequality, that is,
\begin{align*}
\|\Psi(t,\cdot)-\widetilde{T}_{\intt}(|D|)\Phi(t,\cdot)\|_{(\dot{H}^s_q)^3}
&\lesssim\left\|\chi_{\intt}(\xi)|\xi|^s\big(\widehat{\Psi}(t,\xi)-\widetilde{T}_{\intt}\widehat{\Phi}(t,\xi)\big)\right\|_{(L^{q'})^3}\\
&\quad+\left\|\big(1-\chi_{\intt}(\xi)\big)|\xi|^s\left(|\widehat{\Psi}(t,\xi)|+|\widehat{\Phi}(t,\xi)|\right)\right\|_{(L^{q'})^3},
\end{align*}
 we totally complete the proof of Theorem \ref{Thm-Lp-Lq}.

\section{Singular limits and initial layer}\label{Section-Singular-Limit}\setcounter{equation}{0}
\subsection{Formal expansions with small thermal relaxation}
\hspace{5mm}Our main purpose here is devoted to analyze some influence of small thermal relaxation $0<\tau\ll1$ on asymptotic behavior of the solution to the viscous MGT equation \eqref{Eq_MGT}. Strongly motivated by the Prandtl boundary layer theory and the multi-scale analysis (see, for example, \cite{Wang-Wang-Xin=2010}), we look for an approximate solution to the Cauchy problem \eqref{Eq_MGT} by using the WKB expansions. To be specific, we introduce a formal power series expansion in $\tau$ for the solution
\begin{align}\label{Rep}
	\psi(t,x)=\sum\limits_{j\geqslant 0}\tau^j\left(\psi^{I,j}(t,x)+\psi^{L,j}(\tfrac{t}{\tau},x)\right),
\end{align}
where all of the above terms in the summation are assumed to be smooth. Here, $\psi^{I,j}=\psi^{I,j}(t,x)$ stands for the dominant profile for the $j$-th order expansion and $\psi^{L,j}=\psi^{L,j}(z,x)$ with the new ansatz $z:=\frac{t}{\tau}\in[0,+\infty)$ denotes the other profile decaying to zero as $z\to+\infty$. We shall determine these profiles formally soon afterwards.

First of all, we plug the formal expansion \eqref{Rep} into the equation of the Cauchy problem \eqref{Eq_MGT} to arrive at
\begin{align*}
	0&=\sum\limits_{j\geqslant 0}\tau^{j+1}\left(\psi^{I,j}_{ttt}+\tau^{-3}\psi^{L,j}_{zzz}\right)+\sum\limits_{j\geqslant0}\tau^j\left(\psi^{I,j}_{tt}+\tau^{-2}\psi^{L,j}_{zz}\right)-\sum\limits_{j\geqslant0}\tau^j\left(\Delta\psi^{I,j}+\Delta\psi^{L,j}\right)\notag\\
	&\quad-\delta\sum\limits_{j\geqslant0}\tau^j\left(\Delta\psi_t^{I,j}+\tau^{-1}\Delta\psi_z^{L,j}\right)-\sum\limits_{j\geqslant0}\tau^{j+1}\left(\Delta\psi_t^{I,j}+\tau^{-1}\Delta\psi_z^{L,j}\right).
\end{align*}
Matching the terms with the size of $\tau^{j+1}$, we claim
\begin{align*}
	0&=\left(\psi_{ttt}^{I,j}+\psi_{zzz}^{L,j+3}\right)+\left(\psi_{tt}^{I,j+1}+\psi_{zz}^{L,j+3}\right)-\left(\Delta\psi^{I,j+1}+\Delta\psi^{L,j+1}\right)\\
	&\quad-\delta\left(\Delta\psi_t^{I,j+1}+\Delta\psi_z^{L,j+2}\right)-\left(\Delta\psi_t^{I,j}+\Delta\psi_z^{L,j+1}\right).
\end{align*}
Naturally, we always impose $\psi^{I,j}\equiv0\equiv \psi^{L,j}$ when $j<0$. Letting $\tau\downarrow0$, i.e. $z\to+\infty$, the profiles $\psi^{L,j}$ and their derivatives tend to zero so that
\begin{align}\label{Eq_0002}
	\psi_{tt}^{I,j+1}-\Delta\psi^{I,j+1}-\delta\Delta\psi_t^{I,j+1}=\Delta\psi_t^{I,j}-\psi_{ttt}^{I,j},
\end{align}
as well as
\begin{align}
	\psi_{zzz}^{L,j+3}+\psi_{zz}^{L,j+3}=\delta\Delta\psi_z^{L,j+2}+\Delta\psi^{L,j+1}+\Delta\psi_z^{L,j+1},\label{Eq_0003}
\end{align}
whose initial conditions will be derived later. We emphasize that the quantities (known functions) on the right-hand sides of the equations \eqref{Eq_0002} and \eqref{Eq_0003}  can be determined by the matches with the sizes of $\tau^j$ and $\tau^{j-1}$. In other words, we may regard these quantities as the source terms.

Moreover, one notices
\begin{align*}
	\psi_0(x)&=\psi^{I,0}(0,x)+\psi^{L,0}(0,x)+\sum\limits_{j\geqslant 1}\tau^j\left(\psi^{I,j}(0,x)+\psi^{L,j}(0,x)\right),\\
	\psi_1(x)&=\tau^{-1}\psi_z^{L,0}(0,x)+\psi_t^{I,0}(0,x)+\psi_z^{L,1}(0,x)+\sum\limits_{j\geqslant 1}\tau^j\left(\psi_t^{I,j}(0,x)+\psi_z^{L,j+1}(0,x)\right),\\
	\psi_2(x)&=\tau^{-2}\psi_{zz}^{L,0}(0,x)+\tau^{-1}\psi_{zz}^{L,1}(0,x)+\psi_{tt}^{I,0}(0,x)+\psi_{zz}^{L,2}(0,x)+\sum\limits_{j\geqslant 1}\tau^j\left(\psi_{tt}^{I,j}(0,x)+\psi_{zz}^{L,j+2}(0,x)\right).
\end{align*}
Due to  the initial conditions of the Cauchy problem \eqref{Eq_MGT} being independent of $\tau$, we then may automatically consider $\psi^{I,j}(0,x)=0= \psi^{L,j}(0,x)$, $\psi_t^{I,j}(0,x)=0=\psi_z^{L,j+1}(0,x)$ and $\psi_{tt}^{I,j}(0,x)=0=\psi_{zz}^{L,j+2}(0,x)$ when $j\geqslant1$. Consequently, it results
\begin{enumerate}[(1)]
	\item $\psi^{I,0}(0,x)+\psi^{L,0}(0,x)=\psi_0(x)$,
	\item $\psi_z^{L,0}(0,x)=0$ and $\psi_t^{I,0}(0,x)+\psi_z^{L,1}(0,x)=\psi_1(x)$,
	\item $\psi_{zz}^{L,0}(0,x)=0=\psi_{zz}^{L,1}(0,x)$ and $\psi_{tt}^{I,0}(0,x)+\psi_{zz}^{L,2}(0,x)=\psi_2(x)$.
\end{enumerate}
In fact, the initial conditions in the last three summations need to be determined by their corresponding Cauchy problems carefully rather than artificial choices due to the incompatibility of the initial data.

\subsection{Formal derivations of small thermal relaxation profiles}
\hspace{5mm}To begin with our deduction, we take $j=-1$ in \eqref{Eq_0002} and $j=-3$ in \eqref{Eq_0003}, respectively, to find
\begin{align*}
	\begin{cases}
		\psi_{tt}^{I,0}-\Delta\psi^{I,0}-\delta\Delta\psi_t^{I,0}=0,&x\in\mb{R}^n,\ t>0,\\
		\psi^{I,0}(0,x)=\psi_0^{I,0}(x),\ \psi_t^{I,0}(0,x)=\psi_1^{I,0}(x),&x\in\mb{R}^n,
	\end{cases}
\end{align*}
moreover,
\begin{align}\label{Eq_0004}
	\begin{cases}
		\psi_{zzz}^{L,0}+\psi_{zz}^{L,0}=0,&x\in\mb{R}^n,\ z>0,\\
		\psi^{L,0}(0,x)=\psi_0^{L,0}(x),\ \psi_z^{L,0}(0,x)=\psi_{zz}^{L,0}(0,x)=0,&x\in\mb{R}^n,
	\end{cases}
\end{align}
carrying the restriction $\psi_0^{I,0}(x)+\psi_0^{L,0}(x)=\psi_0(x)$, where we used $\psi^{I,j}\equiv0\equiv\psi^{L,j}$ when $j<0$. In order to determine the initial conditions in the last two problems by themselves, let us introduce the auxiliary function $\psi^{P,0}=\psi^{P,0}(t,x)$ such that
\begin{align}\label{Expression-psi-P-0}
	\psi^{P,0}(t,x):=\psi^{L,0}(\tfrac{t}{\tau},x)+\psi_0^{I,0}(x),
\end{align}
which solves the following Cauchy problem:
\begin{align*}
	\begin{cases}
		\tau\psi_{ttt}^{P,0}+\psi_{tt}^{P,0}=0,&x\in\mb{R}^n,\ t>0,\\
		\psi^{P,0}(0,x)=\psi_0(x),\ \psi_t^{P,0}(0,x)=\psi_{tt}^{P,0}(0,x)=0,&x\in\mb{R}^n.
	\end{cases}
\end{align*}
As a consequence, the last Cauchy problem is solved uniquely by $\psi^{P,0}(t,x)=\psi_0(x)$. Again with $\tau\downarrow0$, the expression \eqref{Expression-psi-P-0} gives $\psi_0^{I,0}(x)=\psi_0(x)$ as well as $\psi_0^{L,0}(x)=0$, then the Cauchy problem \eqref{Eq_0004} yields $\psi^{L,0}\equiv0$.  Next, with the aid of $j=-2$ in \eqref{Eq_0003}, one observes
\begin{align}\label{Eq_0005}
	\begin{cases}
		\psi_{zzz}^{L,1}+\psi_{zz}^{L,1}=0,&x\in\mb{R}^n,\ z>0,\\
		\psi^{L,1}(0,x)=0,\ \psi_z^{L,1}(0,x)=\psi_1^{L,1}(x),\ \psi_{zz}^{L,1}(0,x)=0,&x\in\mb{R}^n,
	\end{cases}
\end{align}
carrying the restriction $\psi^{I,0}_1(x)+\psi^{L,1}_1(x)=\psi_1(x)$.
Similarly, let us introduce another auxiliary function $\psi^{P,1}=\psi^{P,1}(t,x)$ such that
\begin{align*}
	\psi^{P,1}(t,x)=\tau\psi^{L,1}(\tfrac{t}{\tau},x)+t\psi_1^{I,0}(x).
\end{align*}
We directly obtain its model as follows:
\begin{align*}
	\begin{cases}
		\tau\psi_{ttt}^{P,1}+\psi_{tt}^{P,1}=0,&x\in\mb{R}^n,\ t>0,\\
		\psi^{P,1}(0,x)=0,\ \psi_t^{P,1}(0,x)=\psi_1(x),\ \psi_{tt}^{P,1}(0,x)=0,&x\in\mb{R}^n,
	\end{cases}
\end{align*}
whose solution is $\psi^{P,1}(t,x)=t\psi_1(x)$. Letting $\tau\downarrow0$ in the representation of $\psi^{P,1}$, we gain $\psi_1^{I,0}(x)=\psi_1(x)$ as well as $\psi_1^{L,1}(x)=0$. It means that the solution to the Cauchy problem \eqref{Eq_0005} is trivial $\psi^{L,1}\equiv0$. Furthermore, we derive $\psi^{I,0}\equiv\psi^0$ because of $\psi_k^{I,0}(0,x)=\psi_k(x)$ for $k=0,1$, that is the solution to the strongly damped wave equation \eqref{Eq_VEDW}. 

With the purpose of determining higher order profiles, we may take $j=0$ in \eqref{Eq_0002} and $j=-1$ in \eqref{Eq_0003}, respectively, to derive the Cauchy problem \eqref{Kuznetsov_I_1}, whose right-hand side can be re-expressed by $F_1=F_1(t,x)$ such that
\begin{align*}
	F_1:=\Delta\psi_t^0-\psi_{ttt}^0=-\delta\Delta^2\psi^0-\delta^2\Delta^2\psi^0_t
\end{align*} thanks to \eqref{Eq_VEDW}, moreover,
\begin{align*}
	\begin{cases}
		\psi_{zzz}^{L,2}+\psi_{zz}^{L,2}=0,&x\in\mb{R}^n,\ z>0,\\
		\psi^{L,2}(0,x)=\psi_z^{L,2}(0,x)=0,\  \psi_{zz}^{L,2}(0,x)=\psi_2(x)-\Delta\psi_0(x)-\delta\Delta\psi_1(x),&x\in\mb{R}^n,
	\end{cases}
\end{align*}
since the fact that
\begin{align*}
\psi^{I,0}_{tt}(0,x)=\Delta\psi^{I,0}(0,x)+\delta\Delta\psi_t^{I,0}(0,x)=\Delta\psi_0(x)+\delta\Delta\psi_1(x).
\end{align*}
We are able to get its solution via simple calculations, that is,
\begin{align*}
	\psi^{L,2}(\tfrac{t}{\tau},x)=\left(\frac{t}{\tau}-1+\mathrm{e}^{-\frac{t}{\tau}}\right)\big(\psi_2(x)-\Delta\psi_0(x)-\delta\Delta\psi_1(x)\big).
\end{align*}
Lastly, let us consider $j=1$ in \eqref{Eq_0002}, namely, the inhomogeneous Cauchy problem \eqref{Kuznetsov_I_2} occurs. Analogously, its source term $F_2=F_2(t,x)$ can be represented via
\begin{align*}
F_2:=\Delta\psi_t^{I,1}-\psi_{ttt}^{I,1}=-\delta\Delta^2\psi^{I,1}-\delta^2\Delta^2\psi^{I,1}_t-\delta\Delta F_1-\partial_tF_1,
\end{align*}
in which we employed \eqref{Kuznetsov_I_1}.

Summarizing the last discussions associated with the general form \eqref{Rep},  we conclude the next formal WKB expansion of the solution. The rigorous demonstration of this expansion until $j=2$ in the $L^q$ framework will be shown in the subsequent parts.
\begin{prop}\label{Prop-WKB-formal}
	The solution $\psi$ to the Cauchy problem for the viscous MGT equation \eqref{Eq_MGT} with small thermal relaxation $0<\tau\ll 1$ formally has the next asymptotic expansion:
	\begin{align*}
		\psi(t,x)&=\psi^0(t,x)+\tau\psi^{I,1}(t,x)+\tau^2\psi^{I,2}(t,x)+\left(\tau t-\tau^2+\tau^2\mathrm{e}^{-\frac{t}{\tau}}\right)\big(\psi_2(x)-\Delta\psi_0(x)-\delta\Delta\psi_1(x)\big)\\
		&\quad+\sum\limits_{j\geqslant 3}\tau^j\left(\psi^{I,j}(t,x)+\psi^{L,j}(\tfrac{t}{\tau},x)\right),
	\end{align*}
	where $\psi^0$ is the solution to the strongly damped wave equation \eqref{Eq_VEDW}, $\psi^{I,1}$ and $\psi^{I,2}$ solve to the inhomogeneous strongly damped wave equations \eqref{Kuznetsov_I_1} and \eqref{Kuznetsov_I_2}, respectively. Moreover, the higher order remainder terms $\psi^{I,j}$ as well as $\psi^{L,j}$ when $j\geqslant 3$ satisfy the differential equations \eqref{Eq_0002} and \eqref{Eq_0003} with the vanishing initial data. 
\end{prop}

\subsection{The inhomogeneous viscous MGT equation for the error term}
\hspace{5mm}According to the formal WKB expansion of $\psi$ proposed in Proposition \ref{Prop-WKB-formal}, let us introduce the error term $u=u(t,x)$ such that
\begin{align}\label{u-def}
	u(t,x)&:=\psi(t,x)-\big(\psi^0(t,x)+\tau\psi^{I,1}(t,x)+\tau^2\psi^{I,2}(t,x)\big)\notag\\
	&\ \quad-\left(\tau t-\tau^2+\tau^2\mathrm{e}^{-\frac{t}{\tau}}\right)\big(\psi_2(x)-\Delta\psi_0(x)-\delta\Delta\psi_1(x)\big).
\end{align}
Then, it fulfills the following inhomogeneous viscous MGT equation:
\begin{align}\label{Error_Equations}
	\begin{cases}
		\tau u_{ttt}+u_{tt}-\Delta u-(\delta+\tau)\Delta u_t=F^{\tau},&x\in\mb{R}^n,\ t>0,\\
		u(0,x)=u_t(0,x)=0,\ u_{tt}(0,x)=u_2^{\tau}(x),&x\in\mb{R}^n,
	\end{cases}
\end{align}
with the source term $F^{\tau}=F^{\tau}(t,x)$ and the initial data $u_2^{\tau}=u_2^{\tau}(x)$ to be calculated later (they depend on the thermal relaxation $\tau$), where we employed the initial conditions in \eqref{Kuznetsov_I_1} and \eqref{Kuznetsov_I_2} to get $u(0,x)=u_t(0,x)=0$. 

To determine $u^{\tau}_2$, by applying repeatedly the equations of  the Cauchy problems \eqref{Kuznetsov_I_1} as well as \eqref{Kuznetsov_I_2}, we know
\begin{align*}
	u_{tt}
	&=\psi_{tt}-\psi_{tt}^0-\mathrm{e}^{-\frac{t}{\tau}}(\psi_2-\Delta\psi_0-\delta\Delta\psi_1)-\tau\left(\psi_{tt}^{I,1}+\tau\psi_{tt}^{I,2} \right),
\end{align*}
and
\begin{align*}
	\psi_{tt}^{I,1}+\tau\psi_{tt}^{I,2}
	&=\tau(I+\delta\partial_t)\Delta\psi^{I,2}+(I-\tau\delta\Delta)(I+\delta\partial_t)\Delta\psi^{I,1}+\delta(2\tau\delta\Delta-I)\Delta^2\psi^0\\
	&\quad+\delta(\tau-\delta+2\tau\delta^2\Delta)\Delta^2\psi_t^0.
\end{align*}
Thus, the third initial condition is given by
\begin{align*}
	u_2^{\tau}&=-\tau\delta(2\tau\delta\Delta-I)\Delta^2\psi_0-\tau\delta(\tau-\delta+2\tau\delta^2\Delta)\Delta^2\psi_1\\
	&=-\tau^2\delta\Delta^2\left[2\delta\Delta\psi_0+(I+2\delta^2\Delta)\psi_1\right]\\
	&=:\tau^2u_2
\end{align*}
under our technical restriction on the initial data $\Delta^2(\psi_0+\delta\psi_1)=0$, in which the function $u_2=u_2(x)$ is independent of $\tau$ obviously.

To compute the source term $F^{\tau}$, let us recall the hyperbolic operator $\ml{L}_{\mathrm{MGT}}$ for the viscous MGT equation \eqref{Eq_MGT}, which can be rewritten in the next two forms:
\begin{align*}
	\ml{L}_{\mathrm{MGT}}:=\tau\partial_t(\partial_t^2-\Delta)+(\partial_t^2-\Delta-\delta\Delta\partial_t)=(\tau\partial_t^3+\partial_t^2)-(I+(\delta+\tau)\partial_t)\Delta.
\end{align*}
According to the equations of the Cauchy problems \eqref{Eq_MGT}, \eqref{Eq_VEDW}, \eqref{Kuznetsov_I_1} and \eqref{Kuznetsov_I_2}, we immediately observe some cancellations leading to
\begin{align*}
	F^{\tau}&=\ml{L}_{\mathrm{MGT}}\psi-\ml{L}_{\mathrm{MGT}}\psi^0-\tau\ml{L}_{\mathrm{MGT}}\psi^{I,1}-\tau^2\ml{L}_{\mathrm{MGT}}\psi^{I,2}\\
	&\quad-\ml{L}_{\mathrm{MGT}}\left[\left(\tau t-\tau^2+\tau^2\mathrm{e}^{-\frac{t}{\tau}}\right)(\psi_2-\Delta\psi_0-\delta\Delta\psi_1)\right]\\
	&=-\tau^3\partial_t(\partial_t^2-\Delta)\psi^{I,2}+\tau\left(t+\delta-\delta\mathrm{e}^{-\frac{t}{\tau}}\right)\Delta(\psi_2-\Delta\psi_0-\delta\Delta\psi_1),
\end{align*}
where we employed $(\tau\partial_t^3+\partial_t^2)(\tau t-\tau^2+\tau^2\mathrm{e}^{-\frac{t}{\tau}})=0$. By direct computations, we may reduce the higher order time-derivatives into the lower order of them as follows:
\begin{align*}
	\partial_t(\partial_t^2-\Delta)\psi^{I,2}
	&=\delta(I+\delta\partial_t)\Delta^2\psi^{I,2}-\delta(2\delta\Delta+\partial_t+2\delta^2\Delta\partial_t)\Delta^2\psi^{I,1}\\
	&\quad+\delta(I+5\delta^2\Delta+4\delta\partial_t+5\delta^3\Delta\partial_t)\Delta^3\psi^0.
\end{align*}
Hence, from another assumption $\Delta\psi_2=0$, which implies $\Delta(\psi_2-\Delta\psi_0-\delta\Delta\psi_1)=0$, the source term can be explicitly represented  by
\begin{align*}
	F^{\tau}&=\tau^3\delta\left[-(I+\delta\partial_t)\Delta^2\psi^{I,2}+(2\delta\Delta+\partial_t+2\delta^2\Delta\partial_t)\Delta^2\psi^{I,1}\right.\\
	&\quad\qquad\,\left.-(I+5\delta^2\Delta+4\delta\partial_t+5\delta^3\Delta\partial_t)\Delta^3\psi^0\right]\\
	&=:\tau^3F,
\end{align*}
in which the function $F=F(t,x)$ is independent of $\tau$.

 At this moment, we have completely determined the source term $F^{\tau}$ as well as the initial data $u^{\tau}_2$ in the Cauchy problem \eqref{Error_Equations}. Let us next propose an energy estimate uniformly in $\tau$. We apply the partial Fourier transform to the inhomogeneous viscous MGT equation \eqref{Error_Equations} to get the Cauchy problem as follows:
 \begin{align}\label{inhomo_Fourier}
 	\begin{cases}
 		\tau\hat{u}_{ttt}+\hat{u}_{tt}+|\xi|^2\hat{u}+(\delta+\tau)|\xi|^2\hat{u}_t=\tau^3\widehat{F},&\xi\in\mb{R}^n,\ t>0,\\
 		\hat{u}(0,\xi)=\hat{u}_t(0,\xi)=0,\ \hat{u}_{tt}(0,\xi)=\tau^2\hat{u}_2(\xi),&\xi\in\mb{R}^n,
 	\end{cases}
 \end{align}
 carrying the source term $\widehat{F}=\widehat{F}(t,\xi)$ such that
 \begin{align*}
 	\widehat{F}&=\delta\left[-(I+\delta\partial_t)|\xi|^4\widehat{\psi}^{I,2}-(2\delta|\xi|^2-\partial_t+2\delta^2|\xi|^2\partial_t)|\xi|^4\widehat{\psi}^{I,1}\right.\\
 	&\qquad\,\left.+(I-5\delta^2|\xi|^2+4\delta\partial_t-5\delta^3|\xi|^2\partial_t)|\xi|^6\widehat{\psi}^0\right]
 \end{align*}
 and the initial data $\hat{u}_2=\hat{u}_2(\xi)$ such that
 \begin{align}\label{data-u2}
 	\hat{u}_2=\delta|\xi|^4\left[2\delta|\xi|^2\widehat{\psi}_0-(I-2\delta^2|\xi|^2)\widehat{\psi}_1\right].
 \end{align}
 \begin{prop}\label{Pointwise-uniformly-tau}
 	The solution $\hat{u}=\hat{u}(t,\xi)$ to the inhomogeneous Cauchy problem \eqref{inhomo_Fourier} fulfills the next pointwise estimate:
 \begin{align*}
 |\hat{u}|^2\leqslant \tau^6\frac{2(\delta+\tau)}{(\delta-\tau)|\xi|^2}\left(|\hat{u}_2|^2+\frac{1+4\tau(\delta-\tau)|\xi|^2}{4(\delta-\tau)|\xi|^2}\int_0^t|\widehat{F}(\eta,\xi)|^2\mathrm{d}\eta\right)
 \end{align*}
with small thermal relaxation $0<\tau\ll 1$.
 \end{prop}
\begin{proof}
By straightforward computations, the next identity holds:
\begin{align*}
	\frac{1}{2}\frac{\mathrm{d}}{\mathrm{d}t}\left(\left|\frac{1}{2}\hat{u}_t+\tau\hat{u}_{tt}\right|^2\right)&=-\frac{\tau}{2}|\hat{u}_{tt}|^2-\frac{\delta+\tau}{2}|\xi|^2|\hat{u}_t|^2-\frac{1}{4}\mathrm{Re}(\hat{u}_{tt}\bar{\hat{u}}_t)-\frac{1}{2}|\xi|^2\mathrm{Re}(\hat{u}\bar{\hat{u}}_t)\\
	&\quad-\tau(\delta+\tau)|\xi|^2\mathrm{Re}(\hat{u}_t\bar{\hat{u}}_{tt})-\tau|\xi|^2\mathrm{Re}(\hat{u}\bar{\hat{u}}_{tt})+\tau^3\mathrm{Re}\left[\widehat{F}\left(\frac{1}{2}\bar{\hat{u}}_t+\tau\bar{\hat{u}}_{tt}\right)\right].
\end{align*}
To compensate the term $\mathrm{Re}(\hat{u}\bar{\hat{u}}_{tt})$, we take
\begin{align*}
	\frac{1}{2}\frac{\mathrm{d}}{\mathrm{d}t}\left(\frac{\tau}{\delta+\tau}|\xi|^2|\hat{u}+(\delta+\tau)\hat{u}_t|^2\right)&=\frac{\tau}{\delta+\tau}|\xi|^2\mathrm{Re}(\hat{u}_t\bar{\hat{u}})+\tau|\xi|^2|\hat{u}_t|^2+\tau|\xi|^2\mathrm{Re}\,(\hat{u}_{tt}\bar{\hat{u}})\\
	&\quad+\tau(\delta+\tau)|\xi|^2\mathrm{Re}(\hat{u}_{tt}\bar{\hat{u}}_t).
\end{align*}
Handling the rest of product terms, we construct
\begin{align*}
	\frac{1}{2}\frac{\mathrm{d}}{\mathrm{d}t}\left(\frac{\delta-\tau}{2(\delta+\tau)}|\xi|^2|\hat{u}|^2+\frac{1}{4}|\hat{u}_t|^2\right)=\frac{\delta-\tau}{2(\delta+\tau)}|\xi|^2\mathrm{Re}(\hat{u}_t\bar{\hat{u}})+\frac{1}{4}\mathrm{Re}(\hat{u}_{tt}\bar{\hat{u}}_t).
\end{align*}
Summing up the last three equalities and employing Cauchy's inequality, one may claim
\begin{align*}
&\frac{1}{2}\frac{\mathrm{d}}{\mathrm{d}t}\left(\left|\frac{1}{2}\hat{u}_t+\tau\hat{u}_{tt}\right|^2+\frac{\tau}{\delta+\tau}|\xi|^2|\hat{u}+(\delta+\tau)\hat{u}_t|^2+\frac{\delta-\tau}{2(\delta+\tau)}|\xi|^2|\hat{u}|^2+\frac{1}{4}|\hat{u}_t|^2\right)\\
&=-\frac{\tau}{2}|\hat{u}_{tt}|^2-\frac{\delta-\tau}{2}|\xi|^2|\hat{u}_t|^2+\tau^3\mathrm{Re}\left[\widehat{F}\left(\frac{1}{2}\bar{\hat{u}}_t+\tau\bar{\hat{u}}_{tt}\right)\right]\\
&\leqslant\tau^6\frac{1+4\tau(\delta-\tau)|\xi|^2}{8(\delta-\tau)|\xi|^2}|\widehat{F}|^2,
\end{align*}
since $0<\tau< \delta$. An integration of the last inequality over $[0,t]$ yields
\begin{align*}
&\left|\frac{1}{2}\hat{u}_t+\tau\hat{u}_{tt}\right|^2+\frac{\tau}{\delta+\tau}|\xi|^2|\hat{u}+(\delta+\tau)\hat{u}_t|^2+\frac{\delta-\tau}{2(\delta+\tau)}|\xi|^2|\hat{u}|^2+\frac{1}{4}|\hat{u}_t|^2\\
&\leqslant \tau^6\left(|\hat{u}_2|^2+\frac{1+4\tau(\delta-\tau)|\xi|^2}{4(\delta-\tau)|\xi|^2}\int_0^t|\widehat{F}(\eta,\xi)|^2\mathrm{d}\eta\right).
\end{align*}
Reserving our desire term, the proof is complete.
\end{proof}
\begin{remark}\label{Rem-viscosity}
	The condition $0<\tau<\delta$ when $0<\tau\ll 1$ naturally holds. From the physical point of view, the diffusivity of sound can be explained by
	\begin{align*}
		\delta=\kappa\left(\frac{1}{c_{\mathrm{V}}}-\frac{1}{c_{\mathrm{P}}}\right)+\frac{4}{3}\mu_{\mathrm{V}}+\mu_{\mathrm{B}},
	\end{align*}
	where we denote the constant pressure and constant volume by $c_{\mathrm{P}}$ and $c_{\mathrm{V}}$, respectively. Moreover, $\mu_{\mathrm{B}}$ and $\mu_{\mathrm{V}}$ stand for the bulk and shear viscosities. Consequently, for a given fluid in the consideration of acoustic waves, we may take $0<\tau\ll 1$ such that $0<\tau<\delta$ always holds.
\end{remark}
\begin{remark}
From the last proof, one may get
\begin{align*}
|\xi|^2|\hat{u}|^2+|\hat{u}_t|^2\leqslant 4\tau^6\left(|\hat{u}_2|^2+\frac{1+4\tau(\delta-\tau)|\xi|^2}{4(\delta-\tau)|\xi|^2}\int_0^t|\widehat{F}(\eta,\xi)|^2\mathrm{d}\eta\right).
\end{align*}
Thus, the singular limits result in Theorem \ref{Theorem-Singular-Limit} also holds for the time-derivative and the spatial-derivative of the solution $\psi$, whose initial data needs higher regularity $H^{10+M_{p,q,n}}_p\times H^{8+M_{p,q,n}}_p$ due to the additional factor $|\xi|^2$ in the previous pointwise estimate.
\end{remark}
 
\subsection{Some estimates for the source term in the Fourier space}
\hspace{5mm}The present part contributes to some estimates for the source term $F$ with respect to the initial data. Indeed, the source term in the Fourier space can be controlled by
\begin{align*}
|\widehat{F}|&\lesssim\left(|\xi|^4|\widehat{\psi}^{I,2}|+|\xi|^4|\widehat{\psi}^{I,2}_t|\right)+\left(|\xi|^6|\widehat{\psi}^{I,1}|+|\xi|^4|\widehat{\psi}^{I,1}_t|+|\xi|^6|\widehat{\psi}^{I,1}_t|\right)\\
&\quad+\left(|\xi|^6|\widehat{\psi}^0|+|\xi|^8|\widehat{\psi}^0|+|\xi|^6|\widehat{\psi}^0_t|++|\xi|^8|\widehat{\psi}^0_t|\right).
\end{align*}
The unexpressed multiplicative constants throughout this subsection are independent of $\tau$, because we have separated the thermal relaxation from the original source term $F^{\tau}$ in the previous deduction. Recalling the zones $\ml{Z}_{\intt}(\varepsilon_0)$, $\ml{Z}_{\bdd}(\varepsilon_0,N_0)$ and $\ml{Z}_{\extt}(N_0)$ introduced in Section \ref{Section-Lp-Lq}, we next will estimate the source term in these zones, respectively. Among them, thanks to $\varepsilon_0\leqslant|\xi|\leqslant N_0$ for $\xi\in\ml{Z}_{\bdd}(\varepsilon_0,N_0)$, we immediately obtain
\begin{align}\label{Est-BDD}
\chi_{\bdd}(\xi)|\widehat{F}|\lesssim \chi_{\bdd}(\xi)\mathrm{e}^{-ct}(|\widehat{\psi}_0|+|\widehat{\psi}_1|),
\end{align}
where we used exponential stabilities for inhomogeneous strongly damped wave equations.

Let us consider small frequencies $\xi\in\ml{Z}_{\intt}(\varepsilon_0)$. Due to the fact that time-derivatives of solutions to strongly damped wave equations bring the additional factor $|\xi|$, we just need to estimate
\begin{align*}
\chi_{\intt}(\xi)|\widehat{F}|\lesssim\chi_{\intt}(\xi)\left(|\xi|^4|\widehat{\psi}^{I,2}|+|\xi|^6|\widehat{\psi}^{I,1}|+|\xi|^4|\widehat{\psi}^{I,1}_t|+|\xi|^6|\widehat{\psi}^0|\right).
\end{align*}
According to the diffusion wave profile for the strongly damped waves (cf. \cite{Ikehata=2014,Ikehata-Onodera=2017,Chen-Ikehata=2023}), the characteristic roots for the Cauchy problem \eqref{Eq_VEDW} are
\begin{align*}
\lambda_{\pm}=\frac{1}{2}\left(-\delta|\xi|^2\pm i\sqrt{4|\xi|^2-\delta^2|\xi|^4}\,\right)=\pm i|\xi|-\frac{\delta}{2}|\xi|^2+O(|\xi|^3)\ \ \mbox{as}\ \ \xi\in\ml{Z}_{\intt}(\varepsilon_0),
\end{align*}
and the pointwise estimate holds
\begin{align*}
\chi_{\intt}(\xi)|\widehat{\psi}^0|\lesssim\chi_{\intt}(\xi)\mathrm{e}^{-c|\xi|^2t}\left(|\cos(|\xi|t)|\,|\widehat{\psi}_0|+\frac{|\sin(|\xi|t)|}{|\xi|}|\widehat{\psi}_1|\right),
\end{align*}
due to the kernels $\widehat{K}_j=\widehat{K}_j(t,|\xi|)$ with $j=0,1$ for $\xi\in\ml{Z}_{\intt}(\varepsilon_0)$ behaving
\begin{align*}
\widehat{K}_0:=\frac{\lambda_+\mathrm{e}^{\lambda_-t}-\lambda_-\mathrm{e}^{\lambda_+t}}{\lambda_+-\lambda_-}\sim\cos(|\xi|t)\mathrm{e}^{-\frac{\delta}{2}|\xi|^2t}\ \ \mbox{and}\ \ \widehat{K}_1:=\frac{\mathrm{e}^{\lambda_+t}-\mathrm{e}^{\lambda_-t}}{\lambda_+-\lambda_-}\sim \frac{\sin(|\xi|t)}{|\xi|}\mathrm{e}^{-\frac{\delta}{2}|\xi|^2t}.
\end{align*}
From Duhamel's principle, the solution for the Cauchy problem \eqref{Kuznetsov_I_1} is shown by
\begin{align*}
\widehat{\psi}^{I,1}&=-\int_0^t\widehat{K}_1(t-\eta,|\xi|)\left[\delta|\xi|^4\left(\widehat{K}_0(\eta,|\xi|)\widehat{\psi}_0+\widehat{K}_1(\eta,|\xi|)\widehat{\psi}_1\right)\right.\\
&\qquad\qquad\qquad\qquad\qquad\left.+\delta^2|\xi|^4\left(\partial_t\widehat{K}_0(\eta,|\xi|)\widehat{\psi}_0+\partial_t\widehat{K}_1(\eta,|\xi|)\widehat{\psi}_1\right)\right]\mathrm{d}\eta.
\end{align*} 
The asymptotic behavior of these kernels results
\begin{align*}
\chi_{\intt}(\xi)|\widehat{\psi}^{I,1}|&\lesssim \chi_{\intt}(\xi)|\xi|^4\left|\int_0^t\widehat{K}_1(t-\eta,|\xi|)\widehat{K}_1(\eta,|\xi|)\mathrm{d}\eta\right|(|\widehat{\psi}_0|+|\widehat{\psi}_1|)\\
&\lesssim \chi_{\intt}(\xi)|\xi|^4\left(|\xi|^{-2}t\mathrm{e}^{-c|\xi|^2t}+|\xi|^{-3}\mathrm{e}^{-c|\xi|^2t}\right)(|\widehat{\psi}_0|+|\widehat{\psi}_1|)\\
&\lesssim \chi_{\intt}(\xi)	\mathrm{e}^{-c|\xi|^2t}(|\widehat{\psi}_0|+|\widehat{\psi}_1|),
\end{align*}
where we noticed
\begin{align*}
\int_0^t\widehat{K}_1(t-\eta,|\xi|)\widehat{K}_1(\eta,|\xi|)\mathrm{d}\eta&=t\frac{\mathrm{e}^{\lambda_+t}+\mathrm{e}^{\lambda_-t}}{(\lambda_+-\lambda_-)^2}+2\frac{\mathrm{e}^{\lambda_-t}-\mathrm{e}^{\lambda_+t}}{(\lambda_+-\lambda_-)^3},\\
\int_0^t\widehat{K}_1(t-\eta,|\xi|)\partial_t\widehat{K}_1(\eta,|\xi|)\mathrm{d}\eta&=t\frac{\lambda_+\mathrm{e}^{\lambda_+t}+\lambda_-\mathrm{e}^{\lambda_-t}}{(\lambda_+-\lambda_-)^2}+\frac{(\lambda_++\lambda_-)(\mathrm{e}^{\lambda_-t}-\mathrm{e}^{\lambda_+t})}{(\lambda_+-\lambda_-)^3}.
\end{align*}
We also obtain
\begin{align*}
\chi_{\intt}(\xi)|\widehat{\psi}^{I,1}_t|&\lesssim \chi_{\intt}(\xi)|\xi|^4\left|\int_0^t\widehat{K}_1(t-\eta,|\xi|)\partial_t\widehat{K}_1(\eta,|\xi|)\mathrm{d}\eta\right|(|\widehat{\psi}_0|+|\widehat{\psi}_1|)\\
&\lesssim \chi_{\intt}(\xi)	|\xi|\mathrm{e}^{-c|\xi|^2t}(|\widehat{\psi}_0|+|\widehat{\psi}_1|).
\end{align*}
Analogously, the solution to the Cauchy problem \eqref{Kuznetsov_I_2} is represented by
\begin{align*}
\widehat{\psi}^{I,2}=-\int_0^t\widehat{K}_1(t-\eta,|\xi|)\left(\delta|\xi|^4\widehat{\psi}^{I,1}(\eta,\xi)+\delta^2|\xi|^4\widehat{\psi}_t^{I,1}(\eta,\xi)-\delta|\xi|^2\widehat{F}_1(t,\xi)+\partial_t\widehat{F}_1(\eta,\xi)\right)\mathrm{d}\eta,
\end{align*}
whose dominant term for small frequencies is
\begin{align*}
\chi_{\intt}(\xi)|\widehat{\psi}^{I,2}|&\lesssim\chi_{\intt}(\xi)|\xi|^4\left|\int_0^t\widehat{K}_1(t-\eta,|\xi|)\left(\widehat{\psi}^{I,1}(\eta,\xi)-\widehat{\psi}^0_t(\eta,\xi)\right)\mathrm{d}\eta\right|\\
&\lesssim\chi_{\intt}(\xi)|\xi|^8\left|\int_0^t\widehat{K}_1(t-\eta,|\xi|)\int_0^{\eta}\widehat{K}_1(\eta-s,|\xi|)\widehat{K}_1(s,|\xi|)\mathrm{d}s\mathrm{d}\eta\right|(|\widehat{\psi}_0|+|\widehat{\psi}_1|)\\
&\quad+\chi_{\intt}(\xi)|\xi|^4\left|\int_0^t\widehat{K}_1(t-\eta,|\xi|)\partial_t\widehat{K}_1(\eta,|\xi|)\mathrm{d}\eta\right|(|\widehat{\psi}_0|+|\widehat{\psi}_1|)\\
&\lesssim \chi_{\intt}(\xi)|\xi|\mathrm{e}^{-c|\xi|^2t}(|\widehat{\psi}_0|+|\widehat{\psi}_1|),
\end{align*}
thanks to the equality
\begin{align*}
\int_0^t\widehat{K}_1(t-\eta,|\xi|)\int_0^{\eta}\widehat{K}_1(\eta-s,|\xi|)\widehat{K}_1(s,|\xi|)\mathrm{d}s\mathrm{d}\eta=t^2\frac{\mathrm{e}^{\lambda_+t}-\mathrm{e}^{\lambda_-t}}{2(\lambda_+-\lambda_-)^3}-3t\frac{\mathrm{e}^{\lambda_+t}+\mathrm{e}^{\lambda_-t}}{(\lambda_+-\lambda_-)^4}+6\frac{\mathrm{e}^{\lambda_+t}-\mathrm{e}^{\lambda_-t}}{(\lambda_+-\lambda_-)^5}.
\end{align*}
In conclusion, the source term can be estimated by
\begin{align}\label{Est-INT}
\chi_{\intt}(\xi)|\widehat{F}|\lesssim\chi_{\intt}(\xi)|\xi|^5\mathrm{e}^{-c|\xi|^2t}(|\widehat{\psi}_0|+|\widehat{\psi}_1|).
\end{align}

Let us turn to large frequencies $\xi\in\ml{Z}_{\extt}(N_0)$. From the same reason as the small frequency part, we should estimate the next terms only:
\begin{align*}
\chi_{\extt}(\xi)|\widehat{F}|\lesssim\chi_{\extt}(\xi)\left(|\xi|^4|\widehat{\psi}^{I,2}_t|+|\xi|^6|\widehat{\psi}_t^{I,1}|+|\xi|^8|\widehat{\psi}^0_t|\right).
\end{align*}
The characteristic roots for the Cauchy problem \eqref{Eq_VEDW} are
\begin{align*}
	\lambda_{\pm}=\frac{1}{2}\left(-\delta|\xi|^2\pm \sqrt{\delta^2|\xi|^4-4|\xi|^2}\,\right)=\begin{cases}
	-\delta^{-1}+O(|\xi|^{-2})\\
	-\delta|\xi|^2+O(1)
	\end{cases}\ \ \mbox{as}\ \ \xi\in\ml{Z}_{\extt}(N_0).
\end{align*}
Due to the kernels for $\xi\in\ml{Z}_{\extt}(N_0)$ behaving
\begin{align*}
\widehat{K}_0\sim\mathrm{e}^{-\delta^{-1}t}-\delta^{-2}|\xi|^{-2}\mathrm{e}^{-\delta|\xi|^2t}\ \ \mbox{and}\ \ \widehat{K}_1\sim\delta^{-1}|\xi|^{-2}\left(\mathrm{e}^{-\delta^{-1}t}-\mathrm{e}^{-\delta|\xi|^2t}\right),
\end{align*}
by tedious but straightforward computations similar to the small frequency part, we obtain
\begin{align*}
\chi_{\extt}(\xi)|\widehat{\psi}^{I,k}_t|&\lesssim \chi_{\extt}(\xi)\mathrm{e}^{-ct}\left(|\xi|^{2+2k}|\widehat{\psi}_0|+|\xi|^{2k}|\widehat{\psi}_1|\right)
\end{align*}
for $k=0,1,2$. As a result, it gives
\begin{align}\label{EST-EXT}
\chi_{\extt}(\xi)|\widehat{F}|\lesssim \chi_{\extt}(\xi)\mathrm{e}^{-ct}\left(|\xi|^{10}|\widehat{\psi}_0|+|\xi|^8|\widehat{\psi}_1|\right).
\end{align}

\subsection{Rigorous demonstration of singular limits: Proof of Theorem \ref{Theorem-Singular-Limit}}
\hspace{5mm}Recalling that our aim is to estimate the error term $u$ defined in \eqref{u-def}, from Proposition \ref{Pointwise-uniformly-tau} associated with \eqref{data-u2} and \eqref{Est-INT}, we get
\begin{align*}
\chi_{\intt}(\xi)|\hat{u}|^2&\leqslant C\chi_{\intt}(\xi)\tau^6\left(|\xi|^{-2}|\hat{u}_2|^2+|\xi|^{-4}\int_0^t|\widehat{F}(\eta,\xi)|^2\mathrm{d}\eta\right)\\
&\leqslant C\chi_{\intt}(\xi)\tau^6\left[|\xi|^{10}|\widehat{\psi}_0|^2+|\xi|^6|\widehat{\psi}_1|^2+|\xi|^6\int_0^t\mathrm{e}^{-2c|\xi|^2\eta}\mathrm{d}\eta(|\widehat{\psi}_0|^2+|\widehat{\psi}_1|^2)\right]\\
&\leqslant C\chi_{\intt}(\xi)(|\xi|^4|\widehat{\psi}_0|^2+|\xi|^4|\widehat{\psi}_1|^2).
\end{align*} Hereafter, the constant $C$ is independent of $\tau$, which may be changed from line to line.
By using the Hausdorff-Young inequality and H\"older's inequality (these steps are the same as those in Subsection \ref{Sub-Est}), it implies
\begin{align*}
\|\chi_{\intt}(D)u(t,\cdot)\|_{L^q}&\leqslant\begin{cases}
C\tau^3\|\chi_{\intt}(\xi)|\xi|^2\|_{L^{\frac{p'q'}{p'-q'}}}\|(\psi_0,\psi_1)\|_{L^p\times L^p}&\mbox{when}\ \ p\neq q\\
C\tau^3\|\chi_{\intt}(\xi)|\xi|^2\|_{L^{\infty}}\|(\psi_0,\psi_1)\|_{L^p\times L^p}&\mbox{when}\ \ p= q
\end{cases} \\
&\leqslant C\tau^3\|(\psi_0,\psi_1)\|_{L^p\times L^p}
\end{align*}
with $1\leqslant p\leqslant 2\leqslant q\leqslant+\infty$.

For another, combining with \eqref{Est-BDD} and \eqref{EST-EXT} in the pointwise estimate stated in Proposition \ref{Pointwise-uniformly-tau}, the next estimate holds:
\begin{align*}
\big(1-\chi_{\intt}(\xi)\big)|\hat{u}|
\leqslant C\big(1-\chi_{\intt}(\xi)\big)\tau^3\left(\langle\xi\rangle^{9}|\widehat{\psi}_0|+\langle\xi\rangle^{7}|\widehat{\psi}_1|\right).
\end{align*}
Hence, following the deduction of \eqref{Large-Infin}, we claim
\begin{align*}
&\left\|\big(1-\chi_{\intt}(D)\big)u(t,\cdot)\right\|_{L^q}\leqslant C\left\|\big(1-\chi_{\intt}(\xi)\big)\hat{u}(t,\xi)\right\|_{L^{q'}}\\
&\leqslant C\begin{cases}
	\tau^3\left\|\big(1-\chi_{\intt}(\xi)\big)\langle\xi\rangle^{-M_{p,q,n}}\right\|_{L^{\frac{p'q'}{p'-q'}}}\left(\|\langle\xi\rangle^{9+M_{p,q,n}}\widehat{\psi}_0\|_{L^{p'}}+\|\langle\xi\rangle^{7+M_{p,q,n}}\widehat{\psi}_1\|_{L^{p'}}\right)&\mbox{when}\ \ p\neq q\\
	\tau^3\left(\|\langle\xi\rangle^{9}\widehat{\psi}_0\|_{L^{p'}}+\|\langle\xi\rangle^{7}\widehat{\psi}_1\|_{L^{p'}}\right)&\mbox{when}\ \ p= q
\end{cases}\\
&\leqslant C\tau^3\|(\psi_0,\psi_1)\|_{H^{9+M_{p,q,n}}_p\times H^{7+M_{p,q,n}}_p},
\end{align*}
for any $1\leqslant p\leqslant 2\leqslant q\leqslant+\infty$ and $M_{p,q,n}$ defined in \eqref{Mpq}, which completes the proof of Theorem \ref{Theorem-Singular-Limit}.
\section{Final remarks}\label{Sec-Final-Rem}
\hspace{5mm}Throughout this manuscript, we have investigated some asymptotic profiles for the viscous MGT equation. If one is interested in the inviscid case $\delta=0$, we state some comments for it. The inviscid MGT equation in the whole space $\mb{R}^n$, namely,
\begin{align}\label{Invis-1}
	\begin{cases}
		\tau\varphi_{ttt}+\varphi_{tt}-\Delta\varphi-\tau\Delta\varphi_t=0,&x\in\mb{R}^n,\ t>0,\\
		\varphi(0,x)=\varphi_0(x),\ \varphi_t(0,x)=\varphi_1(x),\ \varphi_{tt}(0,x)=\varphi_2(x),&x\in\mb{R}^n,
	\end{cases}
\end{align}
with $\tau>0$, is easily to be investigated. According to \cite[Section 2]{Chen-Palmieri=2020}, by setting the new variable with $\tau\varphi_t+\varphi$, its solution $\varphi=\varphi(t,x)$ can be explicitly represented by
\begin{align*}
	\varphi(t,x)=\varphi^0(t,x)+\ml{F}^{-1}_{\xi\to x}\left[\frac{\tau^2}{1+\tau^2|\xi|^2}\left(\frac{\sin(|\xi|t)}{\tau|\xi|}-\cos(|\xi|t)+\frac{1}{2}\mathrm{e}^{-\frac{t}{\tau}}\right)\right]\big(\varphi_2(x)-\Delta\varphi_0(x)\big).
\end{align*}
In the above, $\varphi^0=\varphi^0(t,x)$ is the solution to the free wave equation
\begin{align*}
	\begin{cases}
		\varphi^0_{tt}-\Delta\varphi^0=0,&x\in\mb{R}^n,\ t>0,\\
		\varphi^0(0,x)=\varphi_0(x),\ \varphi^0_t(0,x)=\varphi_1(x),&x\in\mb{R}^n,
	\end{cases}
\end{align*}
which is the formal limit model of \eqref{Invis-1} as $\tau=0$. Hereditarily from \eqref{CON-01} with $\delta=0$ for the inviscid case, the identity $\varphi_2(x)=\Delta\varphi_0(x)$ is the threshold condition for the appearance of the initial layer with the fast change factor $\mathrm{e}^{-\frac{t}{\tau}}$ as $\tau\downarrow0$. Particularly, if $\varphi_2(x)=\Delta\varphi_0(x)$, then $\varphi\equiv\varphi^0$ in any sense of convergence. Due to the free wave kernels in the formula of $\varphi$, via the multiplier method in the pseudo-differential zone and the Littman type lemma associated with the Littlewood-Paley decomposition in the hyperbolic zone (cf. \cite[Section 16]{Ebert-Reissig-book}), one may derive $L^p-L^q$ estimates of $\varphi$ and singular limits as $\tau\downarrow0$, nevertheless, this purpose is beyond the scope of this manuscript.
\section*{Acknowledgments}
 Wenhui Chen is supported in part by the National Natural Science Foundation of China (grant No. 12301270, grant No. 12171317), 2024 Basic and Applied Basic Research Topic--Young Doctor Set Sail Project (grant No. 2024A04J0016), Guangdong Basic and Applied Basic Research Foundation (grant No. 2023A1515012044). The author thanks professors Ya-guang Wang (Shanghai Jiao Tong University) and Ryo Ikehata (Hiroshima University) for the enlightenments and valuable suggestions  of the singular layer theory.

\end{document}